\documentclass[a4paper,11pt]{article}
\usepackage[english]{babel}
\usepackage[T1]{fontenc}
\usepackage[latin1]{inputenc}
\usepackage{amsmath}
\usepackage{amssymb}
\usepackage{color}
\usepackage{graphicx}
\usepackage{psfrag}

\newcommand{\xx}{{\mathbf x}}
\newcommand{\vv}{{\mathbf v}}
\newcommand{\uu}{{\mathbf u}}

\newcommand{\df}{f(t,x,v)}
\newcommand{\dfb}{f(t,{\mathbf x},{\mathbf v})}
\newcommand{\xt}{\tilde{x}}

\newcommand{\rd}{\mathbb{R}^d}
\newcommand{\rn}{\mathbb{R}^N}
\newcommand{\den}{\rho(\xx,t)}
\newcommand{\vel}{\uu(\xx,t)}
\newcommand{\tem}{T(\xx,t)}

\newcommand{\ien}{e(\xx,t)}
\newcommand{\lgl}{\langle}
\newcommand{\rgl}{\rangle}
\newcommand{\dt}{\partial t}

\newcommand{\tn}{t^{n}}
\newcommand{\ndt}{n \Delta t}
\newcommand{\dv}{\Delta  v}
\newcommand{\dx}{\Delta x}
\newcommand{\dlt}{\Delta t}
\newcommand{\gam}{\gamma}
\newcommand{\fnij}{f^{n}_{ij}}
\newcommand{\fntij}{\tilde{f}^{n}_{ij}}
\newcommand{\fnpij}{f^{n+1}_{ij}}
\newcommand{\maxij}{M^{(j)} [\tilde{f}^{n}_{i}]}
\newcommand{\maxj}{M_{j}}
\newcommand{\vdom}{\mathcal{K}}
\newcommand{\h}{\mathcal{H}}
\newcommand{\f}{\mathcal{F}}
\newcommand{\lw}{l_{w}^{(k)}}
\newcommand{\pk}{P^{(k)}}
\newcommand{\xp}{x_{i+1/2}}
\newcommand{\xm}{x_{i-1/2}}
\newcommand{\hf}{\hat{f}}
\newcommand{\brf}{\bar{F}}
\numberwithin{equation}{section}
\begin{document}

\title{A New Class of Conservative Large Time Step Methods for the BGK Models of the Boltzmann Equation.}
\author{Pietro Santagati and  Giovanni Russo \\
\small{Department of Mathematics and Computer Science} \\
\small{University of Catania, Italy}}

\date{}
\maketitle
{\abstract {\scriptsize This work is aimed to develop a new class of methods for the BGK model of the
Boltzmann equation. This technique allows to get high order of accuracy both
in space and time, teoretically without CFL stability limitation. It's based on a Lagrangian formulation of the problem:
information is stored on a fixed grid in space and velocity, and the equation is integrated along the characteristics.
The source term is treated implicitly by using a DIRK (Diagonally Implicit Runge Kutta) scheme in order
to avoid the time step restriction due to stiff relaxation. 
In particular some L-stable schemes are tested by smooth and Riemann problems, both in rarefied and 
fully fluid regimes. Numerical results show good accuracy and efficiency of the method.}}   
\vspace{.3cm}
\section{Introduction} \label{S:intro}
The Kinetic Theory of gases is based on the Boltzmann equation (BE) for the distribution function $\df$, which depends on the time $t$, the physical space coordinates vector $\xx$ and the microscopic velocity space coordinates vector $\vv$,and provides an accurate description of rarefied gas flows. The rarefied regime arises once the Knundsen number $Kn=\lambda/L=O(1)$, where $\lambda$ is the mean free path and $L$ is a certain characteristic length. Because of its nonlinearity and multi-dimensionality, it is very difficult to find analytical solutions of BE. In most pratical cases, numerical methods have to be used. The Boltzmann equation has a wide range of applications. 
It is being used to study flows in external atmosphere, in particular flows around spacecrafts. 
This kind of applications are characterized by so large Knundsen numbers that the kinetic approach gets highly demanded.
Recent applications concerns flows in microchannels and Micro Electro Mechanical Systems (MEMS)
\cite{MEMS}. Here the gas works usually at standard conditions, whereas the mean free path lies on submicron scale, like
the characteristic length of the system.

The numerical solution of the Boltzmann equation leads us to be involved in a very challenging topic. This is due 
mainly to the high dimensionality of the 
problem (in the fully 3D case the distribution function depends on seven independent variables), the nonlinear collision term, and the requirement to preserve the collision invariants at a discrete level. Several strategies have been developed 
to tackle the problem, each one being suitable for appropriate circumstances. The numerical methods may be grouped in two main classes.
One is the class of probabilistic methods, like Direct Simulation Monte Carlo (DSMC) method described by Bird \cite{Bird}, or recently by Pareschi and Russo \cite{PR3}. The second one concerns deterministic methods, 
see Buet \cite{Buet}, Ohwada \cite{Ohwada}, Rogier and Schneider \cite{RS}; for a review on both classes 
of method see \cite{PR4}. 
The most demanding numerically part concerns the computations of the collision term. The peculiarity of the deterministic methods lies in the capability of making us able to get very accurate solutions about special cases. 
Belonging to the deterministic class, spectrally accurate methods have been proposed (see, for example, \cite{PaRu00,FMP06,FiRu03}), 
though the computational cost is high.  
The Monte Carlo methods, on the other hand, allow computations of the collision term in very efficient way. They are easy-fitting, 
and can handle physical problems arising from strongly nonequilibrium conditions. 
Due to the statistical noise, accurate solutions can be pursued  just by many average steps. 
Therefore, Monte Carlo methods are not efficient for systems near equilibrium, whereas deterministic methods are pretty demanded. 
An exaustive example may be the computation of gas flows in MEMS, see \cite{FFL}, particularly characterized by low velocities. 
Pareschi and Caflisch \cite{PC} recently proposed an alternative approach, by modifing the DMCS. 
However the trouble concerning the statistical noise is not fully worked out yet.

The computation time can be reduced by considering simplified models of the Boltzmann equation, like the BGK model, introduced by Bhatnagar, Gross and Krook \cite{bgk54}.
BGK model shows good properties. One most of all in particular, the convergence to the Euler equation once the Knundsen number 
approaches to zero. 
In some cases it works well enough also far from the equilibrium, see \cite{GarzoSantos}. This model has been extensively theoretically 
investigated \cite{Perthame}, and many numerical computations  have been carried out in order to validate its properties (\cite{CP91},\cite{YangHuang}). Some interesting applications of BGK model are described in \cite{Aoki90} and \cite{Aoki97}. Moreover BGK model is also applied to particular flows in nanostructure \cite{GuoZaoShi}.
The classical schemes tipically suffer the inefficiency due to the stiffness arising from the relaxation time getting smaller and smaller.
Standard schemes require the solution of nonlinear systems derived by discrete formulations of the problem in implicit form. Recently this problem has been circumvented by Puppo and Pieraccini \cite{Puppo} by an implicit
formulation of the relaxation term that can be explicitly computed. In this work the problem is formulated in semilagrangian form. 
The discrete form is worked out on a fixed square grid into the phase space, whereas the time integration is performed along the characteristics. 
The time integration is implemented by treating implicitly the relaxation term with a technique similar to the one used in \cite{Puppo}. However, because of the Lagrangian nature, the usual CFL stability restriction 
does not apply. In order to perform the time evolution along the characteristics it is needed to reconstruct all components of the distribution function for each space grid node. This task is solved by using high order pointwise WENO reconstruction described in \cite{carlini}. This class of methods join together high order accuracy and high efficiency, because the CFL limitation is definitely avoided. Moreover, they work up to the fluid dynamic regime, though they suffer from accuracy loss when the Knudsen number gets small. 
Although the schemes are not strictly conservative, numerical tests show that the conservation errors are very small for smooth flows. A fully conservative version of the schemes can be constructed. 
This paper is organized as follow: after this introduction, in section 2, we recall the BGK models and some of its properties. Section 3 is devoted to the detailed description of the not conservative schemes, whereas in Section 4 we present the result of some numerical tests. 
In Section 5 we present the conservative version of the method, and show some results pursued by it. Finally, in the last section we give some conclusions. 

\section{The BGK model}
The BGK model, introduced by Bhatnagar, Gross and Krook \cite{bgk54}, is a simplification of the Boltzmann equation where the collisions are modeled by a relaxation of the distribution function \(\df\) towards the Maxwellian.
It consists of the following initial values problem
\begin{equation}
 \begin{aligned} \label{E:bgk}
    \frac{\partial \dfb}{\partial t} + v \cdot \nabla_{x} \dfb &= \frac{1}{\tau} ( M\left[ \dfb \right] - \dfb ), \\
     f(0,\xx,\vv) &= f_{0}(\xx,\vv)\quad t\geq 0,\quad \xx \in \rd, \vv \in \rn.
 \end{aligned}
\end{equation}

Here \(d\geq 1\) and \(N\geq d\) denote the dimensions of the physical and velocity spaces respectively.
\(M[f] = M(v;\{\rho,u,T\})\) is the local Maxwellian computed by the moments of the distribution function \(\dfb\)

\begin{equation} \label{E:max}
   M(v; \{\rho,u,T\}) = \frac{\rho}{(2\pi R T)^{N/2}} \text{exp} \left( -\frac{| \vv-\uu |^2}{2 R T} \right).
\end{equation}

where \( \rho = \den\), \(\uu = \vel\) and \( T = \tem\) denote the macroscopic fields, namely: density, mean velocity and temperature, which are related to the function $f$ as follows.
Let \(\phi(\vv)=(1,\vv,1/2 \vv^2)^T\) denote the vector of the collision invariants of the distribution function $\dfb$. The moments are given by
\begin{equation} \label{Emom0}
   (\rho,\rho \uu, E)^T = \lgl f \phi(\vv) \rgl,
\end{equation}
where 
\begin{equation} \label{E:mom}
   \langle g \rangle = \int_{\rn} g(\vv) \, d\vv, \qquad g:\rn \mapsto \mathbb{R} .
\end{equation}
The quantity \(E(\xx,t)\) is the total energy and it is related to the temperature by the internal energy \(\ien\)
\[
   \ien=\frac{N}{2} R\tem, \quad \rho e = E - \frac{1}{2} \rho \uu^2.
\]

\subsection{Conservation and Entropy Principle}
Conservation laws for the macroscopic fields are regained by (\ref{E:bgk}) upon multiplication by 
$\phi$ and integration in velocity:
\begin{equation}
	\begin{aligned} \label{E:cons}
	\frac{\partial \lgl f \rgl }{\dt} + \nabla_{x} \cdot \lgl \vv f\rgl &= 0, \\
	\frac{\partial \lgl fv \rgl}{\dt} + \nabla_{x} \cdot \lgl \vv \otimes \vv f\rgl &= 0, \\
      \frac{\partial \lgl \frac{1}{2} \vv^{2} f \rgl}{\dt} + \nabla_{\xx} \cdot \lgl \frac{1}{2} \vv^{2} \vv f\rgl &= 0.
	\end{aligned}
\end{equation}
The Entropy Principle, sometimes called H-theorem, holds also for the BGK models, like for the Boltzmann equation \cite{CercCIP}
\begin{equation} \label{E:Entropy}
	\frac{\partial \lgl f \log f \rgl}{\dt} + \nabla_{\xx} \cdot \lgl \vv f \log f \rgl  \leq 0, \quad \forall \dfb > 0.	
\end{equation}
Once the equilibrium has been established the equality holds, since the distribution function is the Maxwellian.
\newline
BGK models are generally implemented by using \(N=3\), which means that the system is a monoatomic gas with three translational degrees of freedom.
When $ d=1 $ the computation time can be reduced by using the approach proposed in \cite{Aoki90}, preserving the properties of the gas.
In such case, one uses the induced cylindrical symmetry in the velocity space, making the problem properly one dimensional in space and two dimensional in velocity.
As in \cite{CP91}, in this work we choose $d=1$, \(N=1\), since our task is to present the 
methods and testing their properties in simple cases. Thus, the system is a gas with one degree of freedom and the integral in velocity space are evaluated  in \(\mathbb{R}\).\newline
The relaxation time of the BGK model is defined by
\begin{equation} \label{E:reltime1}
	\tau^{-1} = c \rho T^{1-\delta};
\end{equation}
in this definition \(\delta\) is the exponent of the viscosity law of the gas (see \cite{ChapCow}). The constant \(c\) is defined by \(c = R T^{\delta}_{ref} / \mu_{ref}\), where \(\mu_{ref}\) is the viscosity at the reference temperature \(T_{ref}\). \newline
In \cite{Aoki90} the authors use expression (\ref{E:reltime1}) with $\delta = 0$, $\tau^{-1} = A \rho$.

For clarity of exposition, in this work we assume that the collision frequency is constant. 
We shall indicate in a remark how to incorporate in the schemes the dependence of the relaxation time on the 
conservative moments.
We reformulate the BGK equation in non-dimensional form, in such a way the relaxation time plays the role of the Knudsen number. 

In the BGK model, with a single relaxation time, the transport coefficients in the fluid regime (\(Kn \ll 1\)) are not correctly predicted. For example, the Prandtl number \(Pr\) is equal to 1, whereas the correct value is \(\frac{2}{3}\). Because of this limitation, several BGK-like models have been proposed (see, for example, Struchtrup \cite{Struch}, Shakov \cite{Shakov}, Liu \cite{Liu}, Hollway \cite{Hol}, Bouchut and Perthame \cite{BP}). For simplicity we will introduce our class of numerical methods using the classical BGK model, remarking that the application to more sofisticated models is possible as well. 

\section{ Description of the method}
\subsection{A basic first order scheme}
The numerical scheme for the solution of Eq. \eqref{E:bgk} is based on the characteristic formulation of the problem \eqref{E:bgk},
\begin{equation}
	\begin{aligned} \label{E:lag}
 	\frac{d \df}{dt} &= \frac{1}{\tau} ( M \left[\df \right] - \df ), \\
	\frac{dx}{dt} &= v,\\
	x(0) &= \xt , \quad f(0,x,v) = f_{0}(\xt,v) \quad t \geq 0, \quad x,v \in \mathbb{R}.
	\end{aligned}
\end{equation}
Here \( x\) becomes a time dependent variable and its equation in \eqref{E:lag} can be integrated immediately. Hence the BGK model may be presented as follow
\begin{equation}
	\begin{aligned} \label{E:lag2}
 	\frac{d \df}{dt} &= \frac{1}{\tau} ( M \left[\df \right] - \df ), \\
	x(t) &= \xt + vt, \\
	\quad f(0,x,v) &= f_{0}(\xt,v) \quad t \geq 0, \quad x,v \in \mathbb{R},
	\end{aligned}
\end{equation}
and \( x\) is given explicitly. 
 For simplicity we assume constant time step \(\Delta t\) and uniform grid in physical space and velocity domain mesh spacing, \(\dx\) and \( \dv\), respectively and denote the grid points by \(\tn = \ndt\), \(x_i=i \dx\), 
 \(i=1,\ldots,N_x\), \( v_j=j \dv\), \(j=-N_v,\ldots,N_v\), where $N_x$ and $2N_v+1$ are the numbers of grid nodes
 in space and velocity, respectively.  
 Let $\fnij$ denote the approximate solution of the problem \eqref{E:lag2} at time $\tn$ in each spatial and velocity node. A first order explicit scheme is given by
\begin{equation}
	\begin{aligned} \label{E:explicit}
	\fnpij &= \fntij + \dlt \h_{ij}^{n}(f), \\
	   x_i &= \xt_i + v_j \dlt , \quad i=1,...,N_x, \quad j=-N_v,\ldots,N_v, \\
	   \text{where} \\
	 \h_{ij}^{n} & = \frac{1}{\tau} (\maxj [\tilde{f}^{n}_i] - \fntij).
	\end{aligned}
\end{equation}
The numerical solution $\fnpij$ at time $(n+1) \dlt$ requires $\tilde{f}^{n}(\xt_i,v_j)$, 
 denoted  by $\fntij$ in eq. \eqref{E:explicit}.
\begin{figure}[htbp]
\centering
\psfrag{tn}{\( \tn\)}
\psfrag{tn1}{\( \tn + \Delta t \)} \psfrag{fn2}{$ f^{(n)}(x_i,v_j) $}
\psfrag{xt}{\(\xt_i\)} \psfrag{x1}{\(x_{i-1}\)} \psfrag{x2}{\(x_{i-2}\)}
\psfrag{fn}{\( \tilde{f}^{(n)}(\xt_i,v_j)\)}
\psfrag{x}{\(x_i\)} \psfrag{fn1}{\( f^{(n+1)}(x_i,v_j)\)}
\psfrag{vgt0}{\( (v_j>0) \)}
\includegraphics[scale=.40,angle=0]{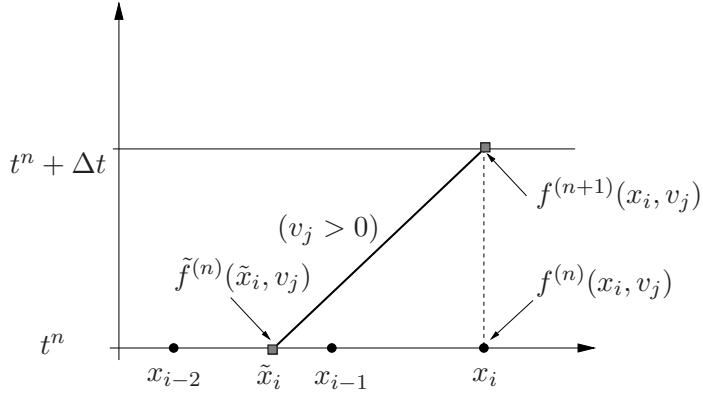}
\caption{Characteristics diagram for positive velocity grid node.} \label{fig:diag1}
\end{figure}
The distribution function's value $\fntij$ is computed by linear interpolation, using the 
values of the distribution function on the left and right nodes of the interval containing $\xt_i=x_i - v_j \dlt$
(see \figurename~\ref{fig:diag1}).

The distribution function for fixed velocity node $v_j$ evolves along the characteristics (thick line). The Maxwellian $\maxij$ is computed as follow
\begin{equation}
	\maxij = \frac{\tilde{\rho}_i}{(2\pi R \tilde{T}_i)^{N/2}} \text{exp} \left( -\frac{ |v-\tilde{u}_i|^2}{2 R \tilde{T}_i} \right).
\end{equation}

This formula requires the computation of the discrete moments of $\fntij$. 
This can be done by using a numerical approximation of the integrals computed in \eqref{E:mom}. Following the notation in \cite{Mieussens}, the discrete velocity grid may be denoted by $\mathcal V$, which is composed of $2N_v+1$ nodes, and $\lgl g \rgl$ can be approximated by a quadrature rule on $\mathcal V$. Let $\lgl g \rgl_{\vdom} $ denote the approximation of $\lgl g \rgl$, where $\vdom$ is the set of $2N_v+1$ indices matching the velocity grid nodes. By this way we compute the moments of the Maxwellian at each grid nodes $\{ x_i \}$,
\[
   (\rho_i,\rho_i u_i, E_i) = \lgl \tilde{f}^{ n }_i \phi(v)\rgl_{\vdom}
\]
As quadrature rule we use summation over $\vdom$, providing spectral accuracy for smooth functions on compact support. The grid $\vdom$ is chosen to include most of the mass. For a given number of nodes $N_v$, an optimal choice of the grid is obtained as a compromise between the extension of the velocity domain and the resolution of the grid. A more sophisticated strategy for the optimal velocity domain configuration lies in adapting the grid to the flow dynamics \cite{RS1},\cite{Santagati}. 

Once the moments are computed on the grid, they can be in turn computed in $\tilde{x}_i$, by a suitable interpolation formula,
so that the Maxwellian gets easily evaluated. Details about the interpolation will be given in Section \ref{sec:WENO}.


The scheme \eqref{E:explicit} can be used to perform the time step. This explicit scheme is first order accurate. Moreover the stability is not preserved as $\tau$ approachs to zero (in this case the sistem \eqref{E:explicit} becomes stiff).
Although the scheme could be made higher order accurate, because of its explicit formulation, it requires very small time steps, 
becoming in turn impractical at the fluid regime. 

In order to circumvent the stiffness arising from the fluid regime, an implicit formulation 
of the system \eqref{E:lag2} is highly desired
\begin{subequations} \label{E:implicit}
	\begin{align}
	 \fnpij &= \fntij + \dlt \h_{ij}^{n+1}(f), \label{E:implicit1}   \\
	   x_i &= \xt_i + v_j \dlt , \quad i=1,...,N_x, \quad j=-N_v,\ldots,N_v, \label{E:implicit2}  \\
	 \h_{ij}^{n+1} &= \frac{1}{\tau} (\maxj [{f}^{n+1}_i] - \fnpij), \label{E:implicit3} .
	\end{align}
\end{subequations}
Hence the moments of $\fnpij$ are needed to compute $\maxj [f^{n+1}_{i}]$, $\forall x_i$ and for each velocity node in $\mathcal V$. As we will see, this can be done explicitly by computing the moments of both sides of \eqref{E:implicit1}
\begin{equation}
  \label{E:fn+1}
  \lgl f^{n+1}(x,v) \phi(v)\rgl = \lgl \tilde{f}^{n}(\xt,v) \phi(v) \rgl + \lgl \mathcal{H}^{n+1}(f) \phi(v)\rgl \dlt.
\end{equation}
As pointed out in \cite{Puppo}, we observe that the moments of the relaxation operator are identically zero 
\[
\lgl \h^{n+1}(f) \phi(v)\rgl = \lgl (M [ {f}^{n+1}(x,v)] - {f}^{n+1}(x,v)) \phi(v) \rgl = 0
\]

therefore

\begin{equation} \label{E:average}
\lgl f^{n+1}(x,v) \phi(v)\rgl = \lgl \tilde{f}^{n}(\xt,v) \phi(v) \rgl, \quad \forall x \in \mathbb{R}.
\end{equation}
and the moments can be easily computed. The Maxwellian $M [ \tilde{f}^{n+1}(x,v)] = 
M(v; \{\rho^{n+1}_i,u^{n+1}_i,T^{n+1}_i\})$ is known and the distribution function at the next time step can be explicitly evaluated

\begin{equation} \label{E:fnp1}
f^{n+1}(x,v) = \frac{\tau \tilde{f}^{n}(x,v) + \dlt M(v; \rho^{n+1}_i,u^{n+1}_i,T^{n+1}_i)}{\tau + \dlt}.
\end{equation}

Once Eqs.~(\ref{E:fn+1}, \ref{E:average}, \ref{E:fnp1}) are discretized on a grid, the resultant first order scheme 
can be written as 

\begin{equation} \label{E:dfnp1}
f^{n+1}_{ij} = \frac{\tau \tilde{f}^{n}_{ij} + \dlt M^{n+1}_{ij}}{\tau + \dlt}.
\end{equation}

where $M^{n+1}_{i,j} = M(v; \rho^{n+1}_i,u^{n+1}_i,T^{n+1}_i$), and the moments are computed by
\begin{eqnarray*}
\rho_{i}^{n+1} & = & \sum_{j=-N_v}^{N_v} \fntij \dv \\
u_{i}^{n+1} & = & \frac{1}{\rho_{i}^{n+1}}\sum_{j=-N_v}^{N_v} v_j \fntij \dv \\
T_{i}^{n+1} & = & \frac{1}{R} \sum_{j=-N_v}^{N_v} v_j^2 \fntij \dv - \left( \rho u^2 \right)_{i}^{n+1} 
\end{eqnarray*}
This approach allows to use CFL numbers greater than one. Moreover, if $\tau \ll 1$ $f^{n}(x,v)$ relaxes very fast toward 
the local Maxwellian.

\noindent {\bf Remark} Since the moments are computed by a quadrature formula, it is not properly true that, in the discrete formulation, $M[\df]$ and $\df$ have the same moments. To get an insight on this aspect see
\cite{Mieussens}. In that paper the author introduces the notion of a discrete Maxwellian, which is more 
consistent with the discrete formulation of the problem. The discrete BGK model obtained using such Maxwellian is 
conservative and entropic. By enough large number of grid points in velocity, the 
continuous and discrete Maxwellians give comparable results. However, for coarse discretization in 
velocity, the discrete Maxwellian introduced in \cite{Mieussens} produces better results.

\medskip

Next section is aimed to present a more general class of methods, based on this basic scheme. The accuracy gets improved, preserving
the discussed properties above.

\subsection{General WENO reconstruction} \label{sec:WENO}
\begin{figure}[htbp]
\centering
\psfrag{xj}{$x_i$} \psfrag{xj1}{$x_{i+1}$}
\psfrag{xjmn}{$x_{i-n+1}$} \psfrag{xjpn}{$x_{i+n}$}
\psfrag{s1}{$S_1$} \psfrag{s2}{$S_n$}
\includegraphics[scale=.40,angle=0]{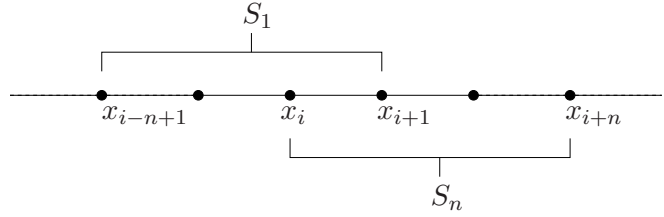}
\caption{Stencil for WENO reconstruction.} \label{fig:diag2}
\end{figure}

The accuracy and the shock capturing properties of the scheme near the fluid regime require a suitable nonlinear reconstruction technique for the computation of $\fntij$.
ENO (Essentially-Non-Oscillstor) and WENO (Weighted-ENO) methods provide the desired high accuracy and non oscillatory properties (see \cite{Shu}). Both methods are based on the reconstruction of piecewise smooth functions by choosing the interpolation points on the smooth side of the function. In ENO methods these points are choosen according to the magnitude of the divided differences evaluated by two candidate stencils. In WENO methods the different polinomials, defined on the stencils, are weighted in such a way that the information about the function on both sides can be used.
Here we focus on WENO reconstruction by introducing the general framework of the implementation. Let us consider the space grid $\{x_i\}_{i \in \mathbb{Z}}$ and the discrete distribution function $\f=\{ f_i\}_{i \in \mathbb{Z}}$ known on any space grid point. For simplicity the time step and velocity grid node indices are not used. The goal is to construct a $2n-1$ degree WENO interpolation on the interval $[x_i,x_{i+1}]$.
Let $L(\xi)$ be the Lagrange polinomial built on the stencil $S=\{x_{i-n+1},...,x_{i+n}\}$. \\
It can be written as follow
\begin{equation} \label{E:lagpol}
L(\xi) = \sum_{k=1}^{n} \lw(\xi) \pk (\xi), \quad \xi \in [x_i,x_{i+1}], \quad k=1,...,n
\end{equation}
where $\lw(\xi)$ are the {\em linear weights}, of $n-1$ degree polinomials, and $\pk(\xi)$
are $n$ degree polinomials, interpolating $\f$ on the stencil $S_k$. 
As shown in \figurename~\ref{fig:diag2}, all the stencils overlap on $[x_i,x_{i+1}]$.
From the definition \eqref{E:lagpol}, the linear weights $\lw$ have to vanish at the nodes outside $S_k$ and must have unit sum in the nodes of $S$, 
where they are nonzero. This condition leads us up to write the linear weights as follow
\begin{equation} \label{E:lw}
	\lw (\xi) = \gamma_{k} \ \prod_{x_l \in S \setminus S_k} (\xi - x_l),
\end{equation}
where $\gamma_{k}$ are evaluated by imposing the unit sum condition
\begin{equation} \label{E:unitsum}
\sum_{k=1}^{n} \lw (x_i) = 1, \quad x_i \in S.
\end{equation}
For a detailed discussion about the calculation of the coefficients $\gamma_k$ we refer to \cite{carlini}.
The WENO reconstruction is expressed by a simple formula
\begin{equation} \label{E:wenopol}
R \left[  \f \right] (\xi) = \sum_{k=1}^{n} w_k(\xi) \pk (\xi).
\end{equation}
Here $R \left[  \f \right] (\xi)$ denotes the numerical reconstruction of the approximate solution, whereas $w_k(\xi)$ are the {\em nonlinear weights}. To ensure stability and consistency the following properties are required
\begin{equation} \label{E:wkcond}
w_k(\xi) \geq 0,\quad  \quad \sum_{k=1}^{n} w_k(x_i) = 1 , \quad i,k=1,...,n.
\end{equation}
Following \cite{Shu} and \cite{JangShu}, such conditions are satisfied by choosing 
\begin{equation} \label{E:nonlinearw}
w_k(\xi) = \frac{\alpha_k(\xi)}{\sum_{l=1}^{n}\alpha_l(\xi)}, \quad k = 1,...,n
\end{equation}
with
\begin{equation} \label{E:alfa}
\alpha_k(\xi) = \frac{\lw(\xi)}{\left( \varepsilon + \beta_k \right) ^2},
\end{equation}
where the parameter $\varepsilon > 0$ is used in order to avoid that the denominator gets 0, whereas $\beta_k$ are the {\em smoothness indicators}.  Usually $\varepsilon=10^{-6}$. Following \cite{Shu}, we use
\begin{equation} \label{E:beta}
\beta_k = \sum_{l=1}^{n} \int_{x_{i}}^{x_{i+1}} \Delta x^{2l-1} \left( \frac{\partial^l \pk (\xi)}{\partial^l \xi} \right)^2 d \xi.
\end{equation}
By this method it is possible to get high order reconstruction in specific points for any cell $[x_i,x_{i+1}]$ in space.

\subsubsection{Third order WENO interpolation}
In this case the stencil is $S=\{x_{i-1},x_i,x_{i+1},x_{i+2}\}$, and the WENO reconstruction operator is given by a superposition of two parabolas $ P^{(1)}(\xi)$ and $P^{(2)}(\xi)$, defined respectively on two overlapping stencils $S_1=\{ x_{i-1},x_i,x_{i+1}\}$ and $S_2=\{x_i,x_{i+1},x_{i+2}\}$,
\begin{equation}
R \left[ \f \right](\xi) = w_1(\xi) P^{(1)}(\xi) + w_2(\xi) P^{(2)} (\xi)
\end{equation}
The linear weights are
\begin{equation}
l_w^{(1)}(\xi) = \frac{x_{i+2}-\xi}{3 \Delta x}, \quad l_w^{(2)}(\xi) = \frac{\xi - x_{i-1}}{3 \Delta x},
\end{equation}
and the smoothness indicators can be computed by \eqref{E:beta}
\begin{equation*}
\begin{aligned} \label{E:beta2}
\beta_1 &= \frac{13}{12} f_{i-1}^2 + \frac{16}{3} f_{i}^2 + \frac{25}{12} f_{i+1}^2 -
	  \frac{13}{3} f_{i-1}f_{i} + \frac{7}{6} f_{i-1}f_{i+1} - \frac{19}{3} f_i f_{i+1} \\
\beta_2 &= \frac{13}{12} f_{i+2}^2 + \frac{16}{3} f_{i+1}^2 + \frac{25}{12} f_{i}^2 -
	   \frac{13}{3} f_{i+2}f_{i+1} + \frac{7}{6} f_{i+2}f_{i} - \frac{19}{3} f_i f_{i+1}.
\end{aligned}	
\end{equation*}
\subsection{Time discretization}
\begin{figure}[!hb]
\centering
\psfrag{tn}{\( \tn\)}
\psfrag{tn1}{\( \tn + \Delta t \)}
\psfrag{xt}{$\xt_i$} \psfrag{Fl2}{$F^{(l-1)}_i$}
\psfrag{tc2}{$t^n + c_{l-1}\dlt$} \psfrag{Fl1}{$F^{(l)}_i$}
\psfrag{x}{\(x_i\)} \psfrag{tc1}{$t^n + c_l \dlt$}
\psfrag{x1}{$x_{i-1}$} \psfrag{x2}{$x_{i-2}$}
\psfrag{vgt0}{\( (v_j>0) \)}
\psfrag{Fc1i}{$\widehat{F}^{(l-1)}_i$}
\psfrag{Fc1j}{$\widehat{F}^{(l-1)}_{i-2}$}
\psfrag{Fc1i1}{$\widehat{F}^{(l-1)}_{i-1}$}
\psfrag{Fc2i}{$\widehat{F}^{(l)}_{i-1}$}
\psfrag{Fc2j}{$\widehat{F}^{(l)}_{i}$}
\includegraphics[scale=.45,angle=0]{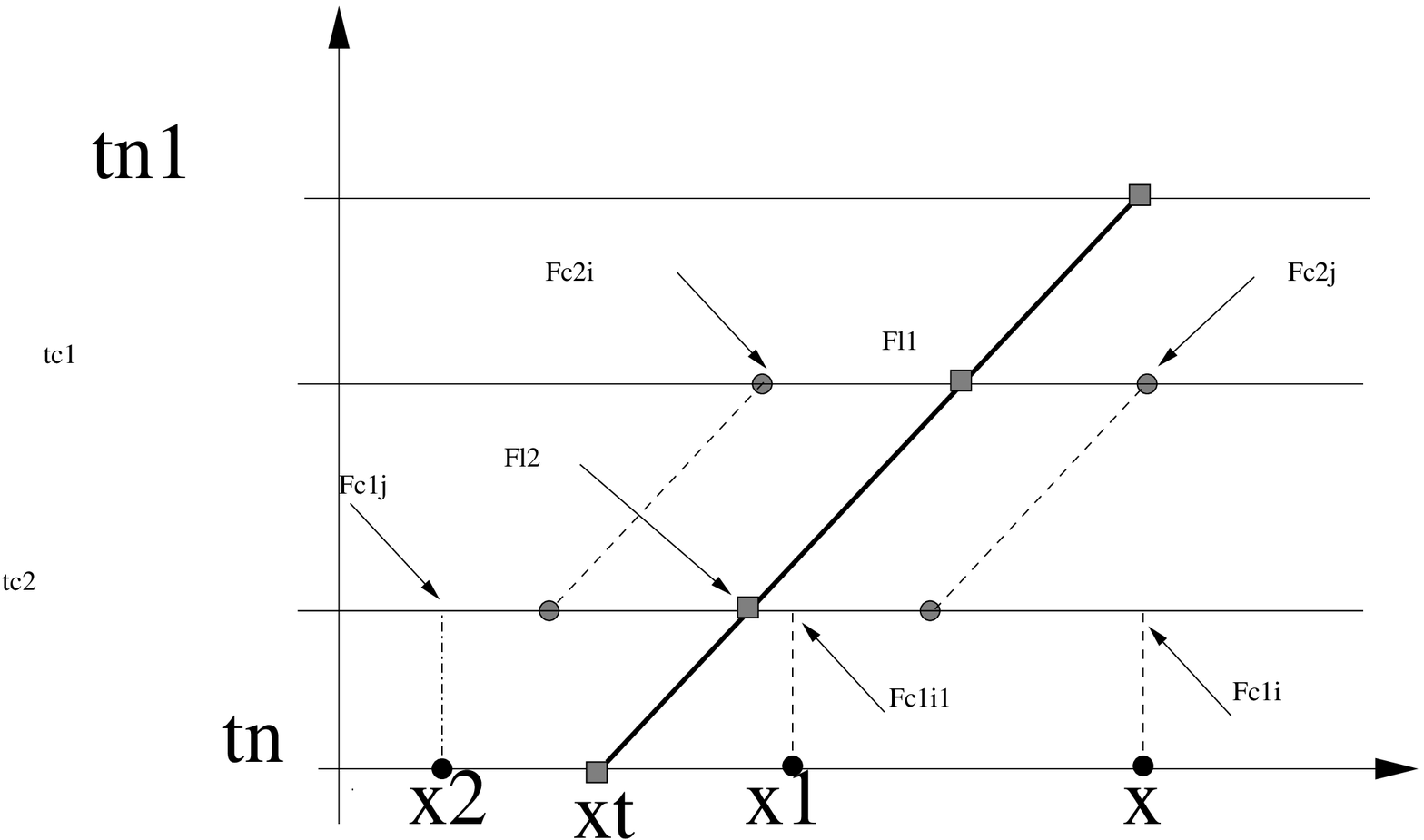}
\caption{Characteristic field for positive velocity grid node.} \label{fig:diag3}
\end{figure}

System \eqref{E:lag2} is a typical ordinary differential equation with relaxation, to be solved in the characteristics framework. Relaxation time lies in a very wide range. It typically extends from order one to very small values compared to the compared to  the  time scale of the  problem. This is the main motivation leading us up to treat the relaxation operator by L-stable diagonally implict Runge Kutta (DIRK) schemes \cite{Hairer,PR1,PR2}. When applied to system \eqref{E:lag2} 
it reads 
\begin{equation}
	\begin{aligned} \label{E:dirk1}
	 F^{(l)}(x,v) &= \tilde{f}^n(\xt,v) + \dlt \sum_{k=1}^{l} a_{lk} \frac{1}{\tau} (M [ F^{(k)}(x,v)] - F^{(k)}(x,v))\\ 
	 f^{n+1}(x,v) &= \tilde{f}^n(\xt,v) + \dlt \sum_{k=1}^{l} w_{k} \frac{1}{\tau} (M [ F^{(k)}(x,v)] - F^{(k)}(x,v)) \\ 
	\end{aligned}
\end{equation}
The triangular $\nu \times \nu$ matrix, $A=(a_{lk})$, and the coefficient vectors, $c=(1,...,c_{\nu})^T$ and $w=(1,...,\nu)^T$, are given by consistency and order conditions. They characterize completely a DIRK scheme, which can be rappresented by the Butcher's {\em tableaux}
\begin{center}
\begin{tabular} { l | c }
$c$ & $A$ \\ \hline
\\
& $w^T$ \\
\end{tabular}
\end{center}
The internal stages are practically evaluated by a sequence of elementary implicit Euler steps.
Scheme \eqref{E:dfnp1} corresponds to \eqref{E:dirk1} with implicit Euler scheme. It will be denoted S1. 
The DIRK methods considered in this work are
\[
    S2 = 
    \begin{array}{c|ccc}
       \alpha &    &  \alpha        &         \\
       $  1 $ &    & 1 - \alpha     &  \alpha \\ \hline
              &    & 1 - \alpha     &  \alpha 
    \end{array}
    \quad \quad
    S3 = 
    \begin{array}{c|cccc}
           1/2 &    &   \gam     \\
    (1+\gam)/2 &    &  (1-\gam)/2  &   \gam    \\
           1   &    & 1 -  \delta  -  \gam  &  \delta  &  \gam  \\ \hline
               &    & 1 -  \delta  -  \gam  &  \delta  &  \gam 
    \end{array}
 \]
which are the implicit parts of the ImEx schemes ARS(2,2,2) and ARS(3,4,3) respectively, see \cite{PR1}, and are 
a second and a third order L-implicit schemes. 
The coefficients are :
\begin{equation*}
 \alpha = 1 - \frac{\sqrt{2}}{2}, \quad \gam=0.4358665215,\quad  \delta=-0.644373171.
\end{equation*}
The schemes constructed by the \emph{tableau} above will be denoted by $S2$ and $S3$ respectively.
The L-stability guarantees that the schemes are able to capture the limit $\tau \to 0$. Other choices are possible of course.
The terms $F^{(l)}(x,v)$ are the {\em internal stages}. We point out that the internal stages have to be evaluated along the characteristics solving an implicit equation, see \figurename~\ref{fig:diag3}. Due to the strong nonlinearity of the relaxation operator, this task, in principle, is not easy to be solved. To circumvent this difficulty we proceed as follow. 
At the starting point all we have is the initial condition. The goal is to evaluate $F^{(1)}_i(v)$, in all discrete velocity domain, on the characteristics framework at the spatial coordinate $x^{(1)}_i = x_i - v(1 - c_1) \dlt$. The prepocessing calculation consists of providing a preliminary internal stage, which is denoted by $\widehat{F}^{(1)}_i(v)$, in each spatial grid node $x_i$, $i = 1,\dots,N_x$. To this end we proceed by performing a single time step of amplitude $c_1\dlt$ using an implicit first order scheme, as shown before.
\begin{figure}[!ht]
\centering
\psfrag{tn}{\( \tn\)}
\psfrag{tn1}{\( \tn + \Delta t \)}
\psfrag{xt}{$\xt_i$}
\psfrag{tc2}{$t^n$}
\psfrag{Fl1}{$F^{(1)}_i$}
\psfrag{x}{\(x_i\)} \psfrag{tc1}{$t^n + c_1 \dlt$}
\psfrag{x1}{$x_{i}-v c_1 \dlt$}
\psfrag{x2}{$x_{i-2}$}
\psfrag{x3}{$x_{i-1}$}
\psfrag{vgt0}{\( (v_j>0) \)}
\psfrag{Fc2i}{$\widehat{F}^{(1)}_{i-1}$}
\psfrag{Fc2j}{$\widehat{F}^{(1)}_{i}$}
\includegraphics[scale=.45,angle=0]{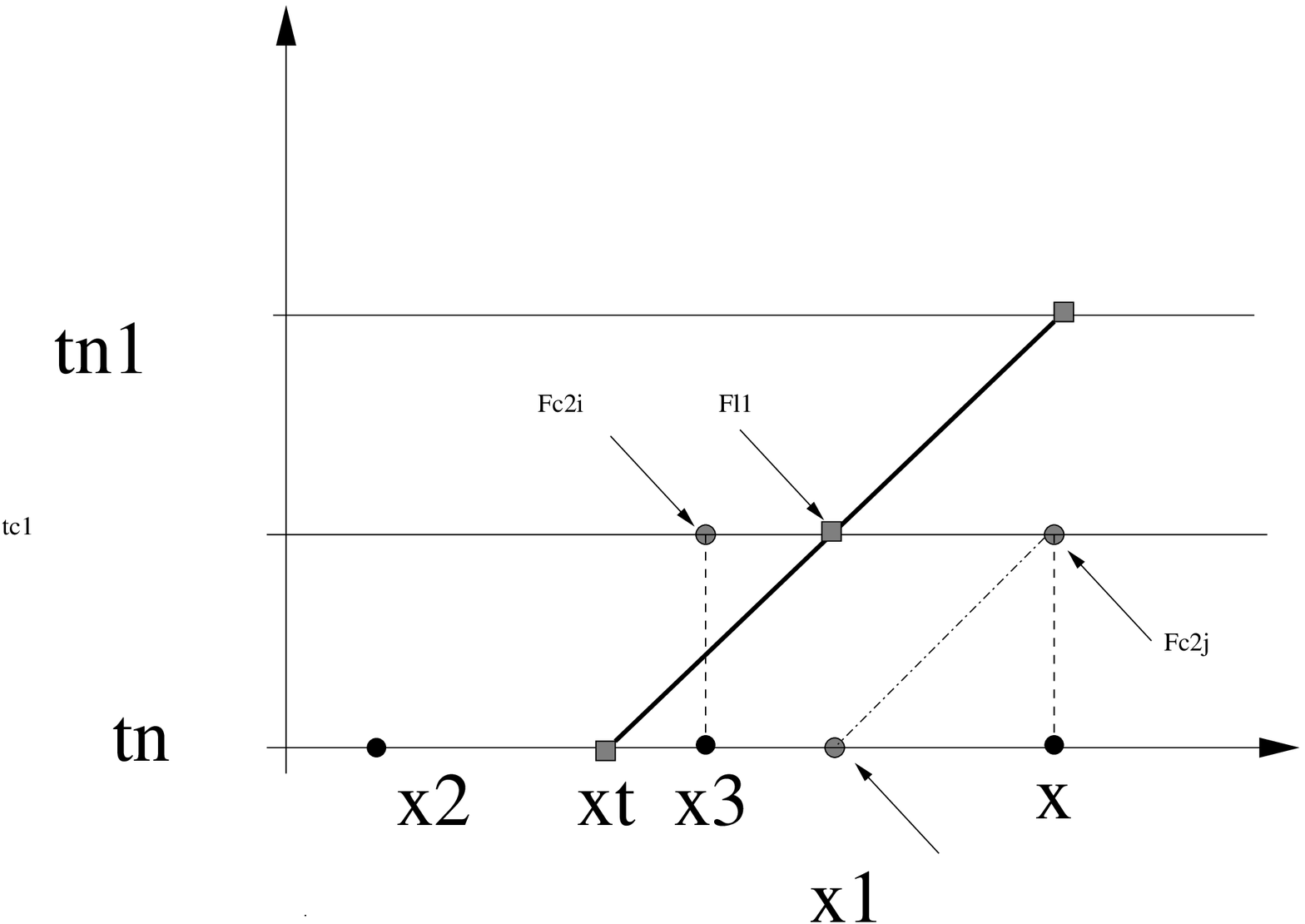}
\caption{Characteristic field for positive velocity grid node.}
\end{figure}
Thus we can write
\begin{equation}
\begin{split} \label{E:firststage}
\widehat{F}^{(1)}_i(v) & = \tilde{f}^n_i(v) +  \frac{\dlt}{\tau} a_{11} (M[\widehat{F}^{(1)}_i](v) - \widehat{F}^{(1)}_i(v)) 
\end{split}
\end{equation}
and \eqref{E:average}
\[
  \lgl \widehat{F}^{(1)}_i(v) \phi(v) \rgl_{\vdom} = \lgl \tilde{f}^n_i(v) \phi(v) \rgl_{\vdom}, \quad i = 1,\ldots,N_x
\]
since
\[
  \lgl (M[\widehat{F}^{(1)}_i](v) - \widehat{F}^{(1)}_i(v)) \phi(v)\rgl_{\vdom} = 0.
\]
Once the moments, $\rho^{(1)},u^{(1)}$ and $T^{(1)}$ are computed, the Maxwellian at the first stage $M[\widehat{F}^{(1)}_i](v),i=1,...,N_x$, is evaluated by
\begin{equation} \label{E:max2}
M[\widehat{F}^{(1)}_i](v) = \frac{\rho^{(1)}_i}{(2\pi R T^{(1)}_i)^{1/2}} \text{exp} \left( -\frac{| v-u^{(1)}_i |^2}{2 R T^{(1)}_i} \right).
\end{equation}
Finally the internal stage $F^{(1)}_i(v)$ is computed by WENO reconstruction at the point $x^{(1)}_i$ by the values $\hat{F}^{(1)}_i(v)$.

For multi-stage case, in order to evaluate the relaxation operator for the $l$-th internal stage
\[
 \frac{1}{\tau}(M [ F^{(l)}_i(v)] - F^{(l)}_i(v))
\]
we need to compute the moments of $\{F^{(l)}_i(v)\},i=1,...,N_x$. To this end, let the set of internal stages
$\{F^{(1)},F^{(2)},...,F^{(l-1)}\}$ be already computed at the characteristic's set points $\{x^{(1)},x^{(2)},\dots,x^{(l-1)}\}$, for each spatial grid node. Therefore we know the pre-processed stages $\{\widehat{F}^{(1)},
\widehat{F}^{(2)},...,\widehat{F}^{(l-1)}\}$. The $l$-th internal stage is evaluated repeating the initial strategy, as showed in \figurename~\ref{fig:diag3}, by using the  pre-processed stages $\widehat{F}^{(l-1)}$ to get the moments
$\rho^{(l)},u^{(l)}$ and $T^{(l)}$ and the Maxwellian $M[\widehat{F}^{(l)}_i](v),i=1,...,N_x$.
Hence it is possible to write
\begin{equation}
      \begin{aligned} \label{E:lstage}
       F^{(l)}_i(v) = f^*_i(v) + a_{ll} \frac{\dlt}{\tau} (M[F^{(l)}_i](v) - F^{(l)}_i(v))
      \end{aligned}
\end{equation}
where
\[
   f^*_i(v) \equiv \tilde{f}^n_i(v) + \dlt \sum_{k=1}^{l-1} a_{lk} (M[F^{(k)}_i](v) - F^{(k)}_i(v))
\]
Equation (\ref{E:lstage}) is then solved for $F^{(l)}_i(v)$ by the same technique used
for Eq.~(\ref{E:firststage}).
Once all internal stages are evaluated, finally we can perform the time evolution step

\noindent {\bf Remark}
In practice the Runge Kutta fluxes can be computed from the internal stages. For example 
\[ 
  \frac{\Delta t}{\tau}(M^{(1)}_{ij}-F^{(1)}_{ij}) = \frac{F^{(1)}_{ij}-\fntij}{a_{11}}.
\]
Hence the scheme can be used in the limit $\tau\to 0$, with no constrain on the time step amplitude. 

\medskip

\begin{equation}
\begin{aligned} \label{E:dirk2}
	 f^{n+1}_i(v) &= \tilde{f}^n_i(v) + \dlt \sum_{k=1}^{\nu} w_{k} \frac{1}{\tau} (M[F^{(k)}_i](v) - F^{(k)}_i(v)).
	\end{aligned}
\end{equation}

\section{Numerical Tests}
These tests are aimed to verify the accuracy (test 1) and the shock capturing properties (test 2) of the schemes.

\subsection{Test 1: regular velocity perturbation}
\begin{figure}[htbp] 
\centering%
 \includegraphics[height=6.5cm,width=6.5cm,angle=0]{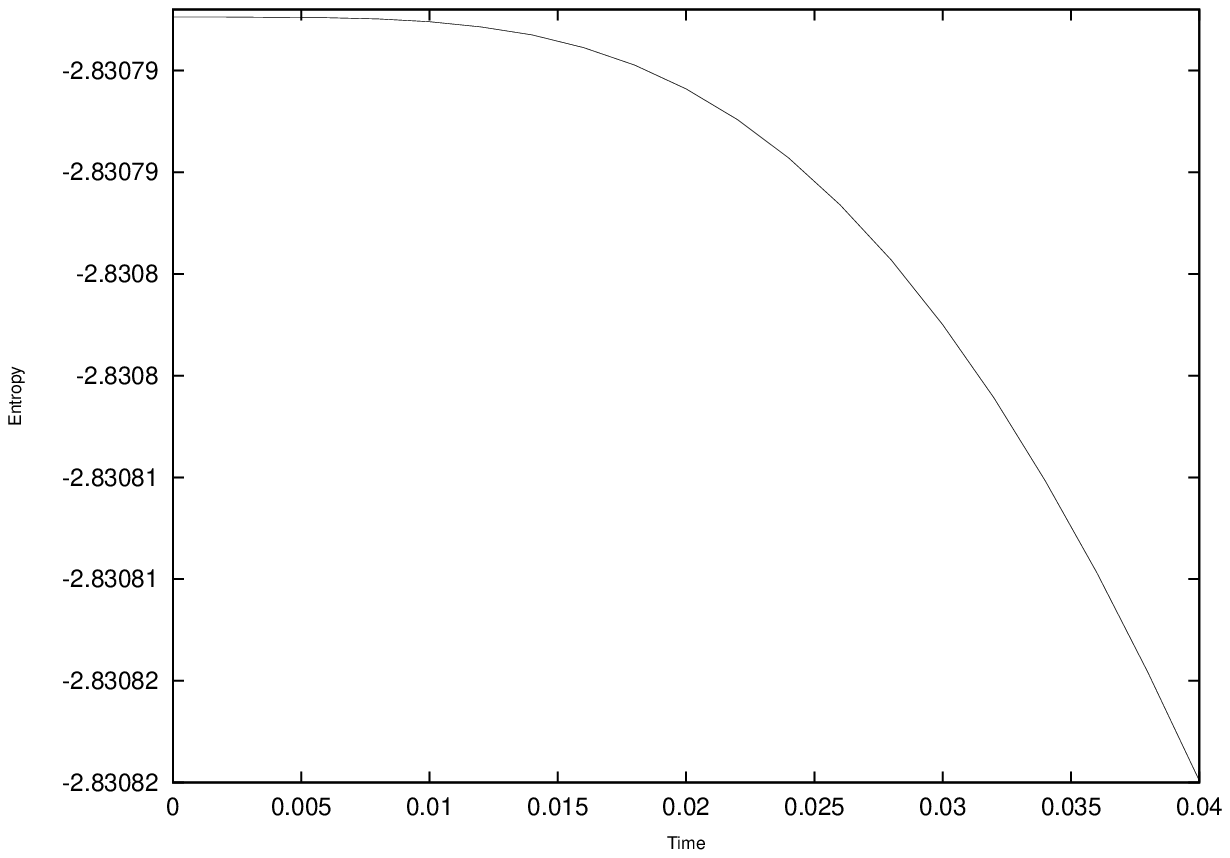}%
 \includegraphics[height=6.5cm,width=6.5cm,angle=0]{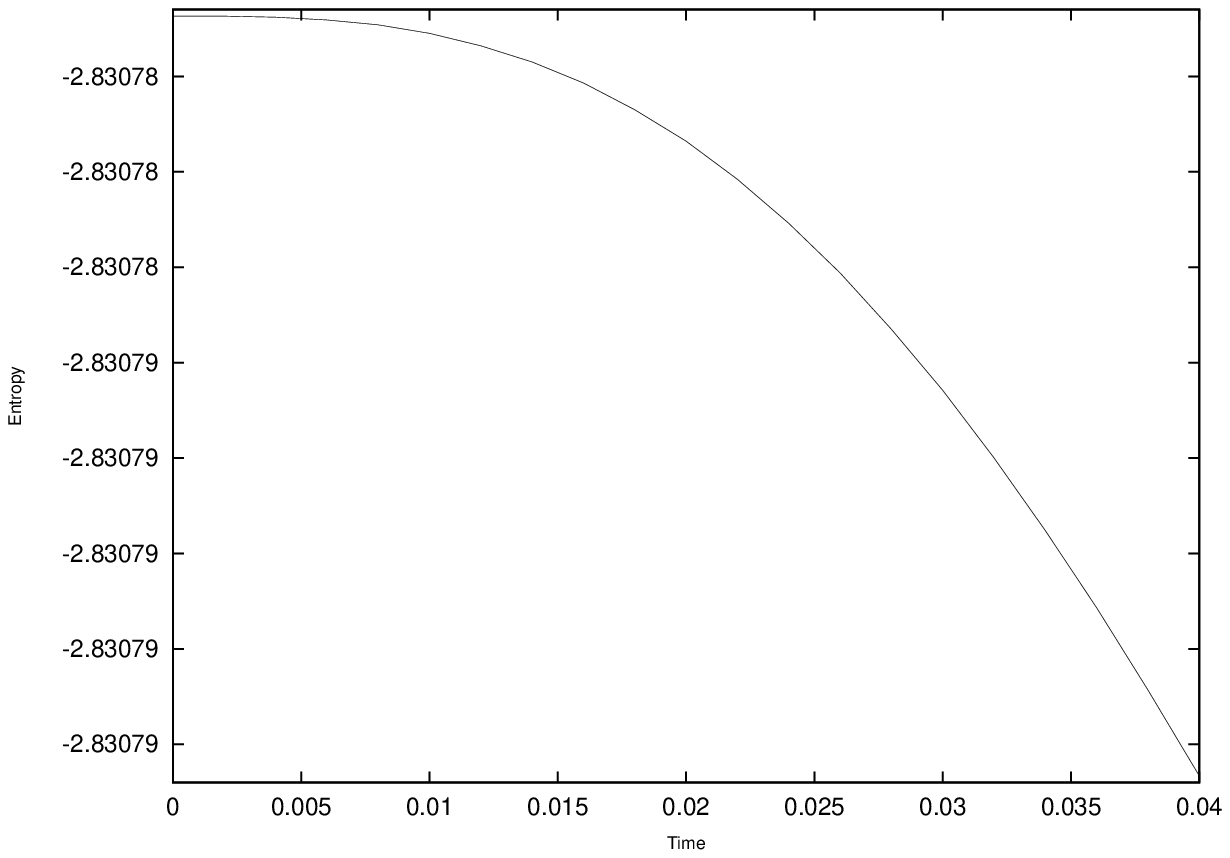}
 \caption{Entropy for test 1;scheme S3, $N_x=800, \tau=10^{-2}$ (left), $\tau= 10^{-6}$ (right).}
 \includegraphics[height=6.5cm,width=6.5cm,angle=0]{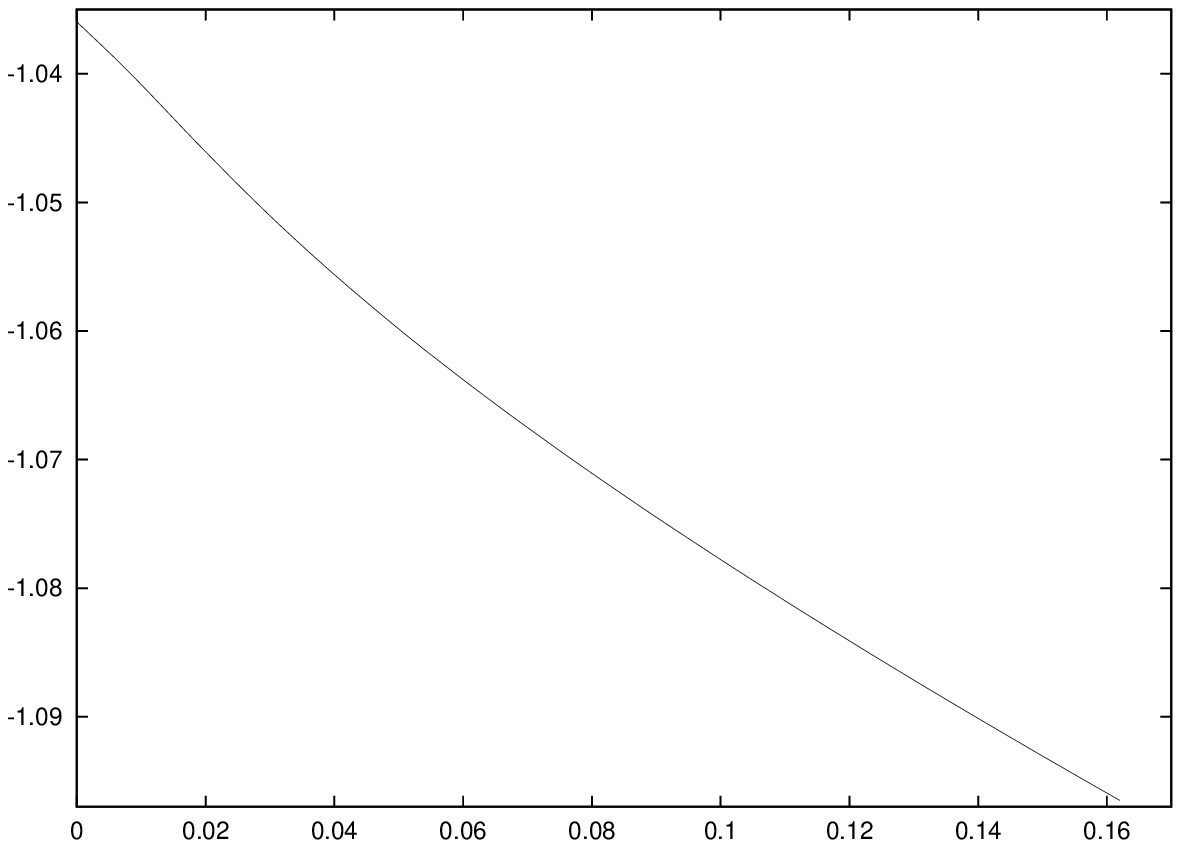}%
 \includegraphics[height=6.5cm,width=6.5cm,angle=0]{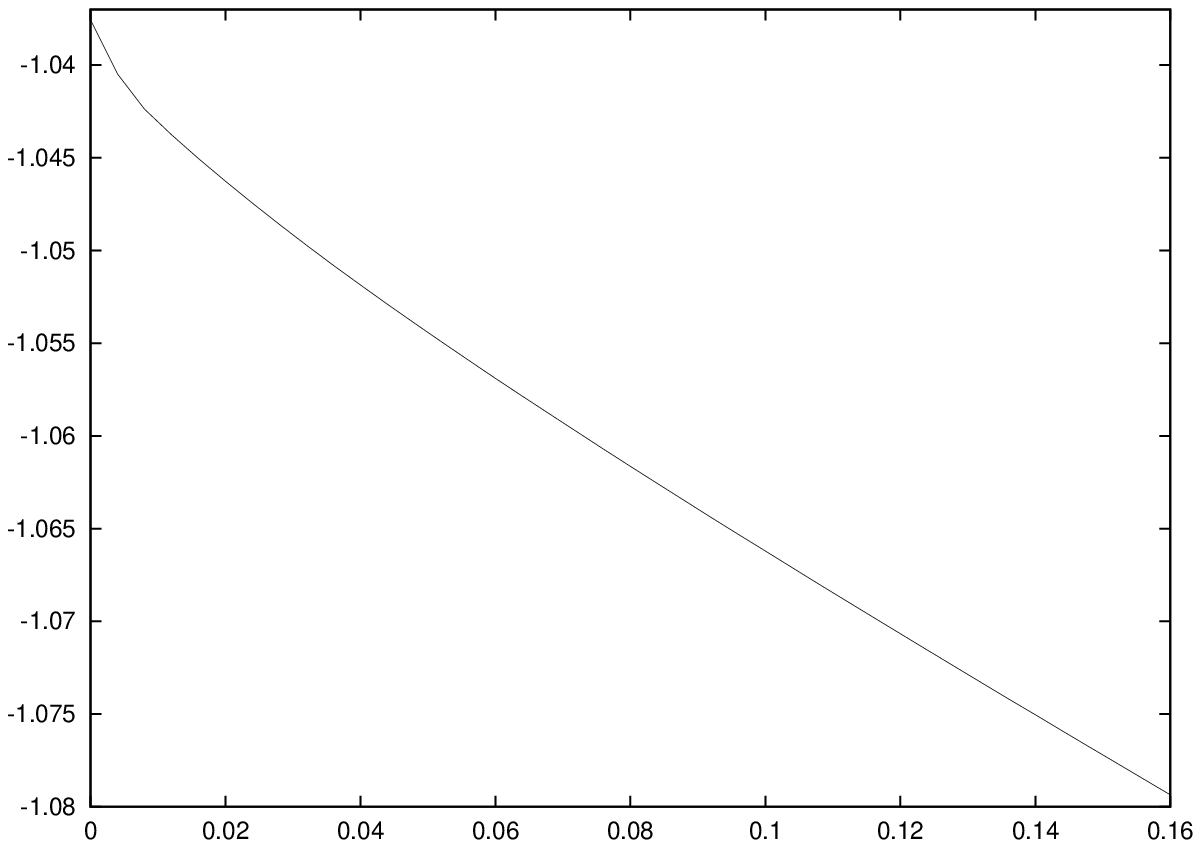}
 \caption{Entropy for test 2; scheme S3, $N_x=800, \tau=10^{-2}$ (left), $\tau= 10^{-6}$ (right).}
 \label{fig:entropy} 
\end{figure}

This test has been proposed in \cite{Puppo}. The solution is smooth, and the accuracy can be tested.
Initial velocity profile is given by
\[
u_0(x) = \frac{1}{\sigma} \left( \exp \left( - (\sigma x -1)^2 \right) - 2 \exp \left( - (\sigma x +3)^2 \right) \right), \quad x \in [-1,1]
\]
where $\sigma$ is a positive constant parameter. Initial density and temperature profiles are uniform, with constant value, $\rho = 1$ and $T=1$ respectively.
The initial condition for the distribution function is the Maxwellian, computed by given macroscopic fields.
The boundary conditions are imposed by prescribed moments as well. 
Two regimes (rarefied and fluid) have been investigated, corresponding to different Knundsen numbers, $\tau=10^{-2}$ and $\tau=10^{-6}$.
The final time, for both cases was 0.04, showed being large enough in order to get the thermodynamic equilibrium.
Accuracy and conservation tests have been performed at the final time. The errors has been computed using a reference solution, defined on a finer grid, with $N_x=1280$ and $N_v=20$. 
The test case has been performed using $N_v=20$ (as for the reference solution),for each spatial grid nodes number, uniformly spaced in [-10,10]. Total entropy 
\begin{equation} \label{E:totentropy}
H(f) = \int_{-1}^{+1} \lgl f \log(f) \rgl \, dx
\end{equation}
has been also computed by a fourth order integration formula in space (see \cite{Puppo}),
for both of the Knundsen numbers and it is reported in \figurename~\ref{fig:entropy}. It is possible to observe that the functional decreases during the time evolution, as expected \cite{CercCIP}. The relative errors and order of accuracy are shown in Tables [\ref{tab:err1},\ref{tab:err2},\ref{tab:err3},\ref{tab:err4}], for the schemes $S2$ and 
$S3$. By using a reference solution at $CFL=0.5$, with $N_x=200$ spatial nodes (whereas $N_v$ was unchanged). Several computation have been carried out, for different CFL numbers. This test is aimed to check the correct behaviour of the schemes as the CFL  changes, leading up to diffusity spurious errors. The results in Table \ref{tab:cflComp} show that the schemes work very good for each order of accuracy, by limiting the errors in a narrow and numerically satisfactory range. Finally the conservation has been investigated. Despite the schemes are not strictly conservative, conservation properties look good, though not many space grid nodes are used. That is because for smooth solutions the weights are close to their linear values. 
A conservative version of this class of schemes is introduced in the next section.

\subsection{Test 2: Riemann problem}
This test allows us to evaluate the capability of our class of schemes in capturing shocks, contact discontinuities and the density profile in a rarefaction. The macroscopic fields are initially assigned in the domain, satisfying the Rankine-Hugoniot shock jump conditions. 
In particular we are interested to the behaviour in the fluid regime limit. The presented results are density, velocity and temperature profiles, for $\tau=10^{-2}$ and $\tau=10^{-6}$, respectively. Like for test 1, the boundary conditions are imposed by Maxwellians computed by prescribed macroscopic moments.
Total entropy calculation results are also presented.
For this test two values $\tau$ are employed , $\tau=10^{-2}$ and $\tau=10^{-6}$. $N_v=60$ nodes are used in the range [-10,10] 
of the discrete velocity domain, as in \cite{Puppo}.

\section{The conservative version}
Here we introduce the conservative version, first order accurate. Generalizing to higher order of accuracy is pretty straightfoward. 
Let us consider the original problem (\ref{E:bgk}). The first internal stage $F^{(1)}(v)$ can be computed by solving the equation \eqref{E:consForm}
\begin{equation} \label{E:consForm}
F^{(1)}_i (v) = f^{n}_i (v) - \dlt \nabla (v F^{(1)}) + \frac{\dlt}{\tau} (M[F^{(1)}_i](v) - F^{(1)}_i(v)).
\end{equation}
Of course a suitable discrete form of the convective part is needed. To overcome this point we use a conservative
finite difference approximation introduced in \cite{ShuOsher}. This efficient method is aimed to look for a function $\hat{f}$ (typically a polynomial) such that the data $\bar{F}^{(1)}=vF^{(1)}$ are interpolated, in the sense of cell average
\[
\bar{F} = \frac{1}{\dx} \int_{\xm}^{\xp} \hat{f} \, dx 
\] 
in such a way that
\begin{equation} \label{shuDer}
\partial _{x} \bar{F}^{(1)}|_{x_i} = \frac{1}{\dx} \left( \hat{f}  (\xp) - \hat{f} (\xm) \right).
\end{equation} 
\begin{figure}[h!]
\centering
\psfrag{tn}{\( \tn\)}
\psfrag{tn1}{\( \tn + \Delta t \)} \psfrag{fn2}{$ f^{(n)}_{i}(v) $}
\psfrag{xt}{\(\xt_i\)} \psfrag{x1}{\(x_{i+\frac{1}{2}}\)} \psfrag{x2}{\(x_{i-\frac{1}{2}}\)}
\psfrag{fn}{$ \tilde{f}^{(n)}_{i}(v) $}
\psfrag{x}{\(x_i\)} \psfrag{fn1}{$ F^{(1)}_{i}(v) $}
\psfrag{vgt0}{\( (v>0) \)}
\includegraphics[scale=.40,angle=0]{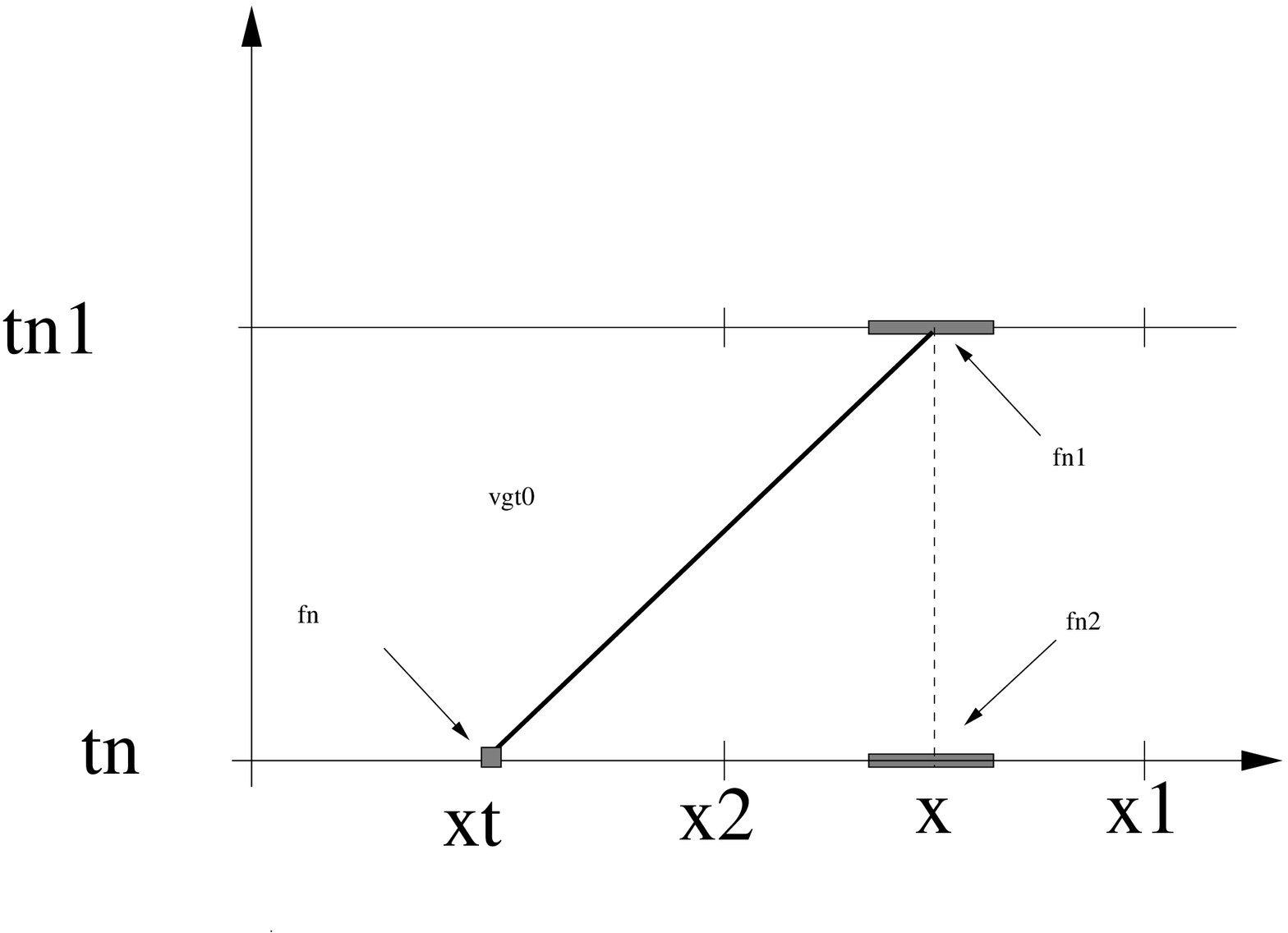}
\caption{Conservative diagram.} \label{fig:diag5}
\end{figure}
We remark that the flux depends linearly on the function $F$ itself. The flux $\hat{f}$ is splitted of course in two contributions 
at each cell border
\begin{equation} \label{E:fluxes}
\hat{f}(\xp) = \hat{f}^{+} (\xp^{-}) + \hat{f}^{-} (\xp^{+}).
\end{equation} 
The reconstruction is performed on $\bar{F}^{+}$ and $\bar{F}^{-}$ that are defined as follow
\begin{equation}
\bar{F}^{+} = 
\begin{cases}
v F^{(1)}, \quad & v>0, \\
0, \quad & v \le 0,
\end{cases} 
\qquad
\bar{F}^{-} = 
\begin{cases}
0, \quad & v \ge 0, \\
v F^{(1)}, \quad & v<0.
\end{cases} 
\end{equation}
Here we use a polynomial reconstruction based on 
WENO method, like in \cite{Puppo} (see \cite{ShuOsher} for the details). \\ Let us consider the stencils 
$\{ \hf_{i+l-1}, \hf_{i+l}, \hf_{i+l+1} \} _{l=-1,0,1}$ and denote by $L^{l}_{i}(x)$ the parabolas that interpolate the given data, 
in the sense of cell average.
Since the flux is linear, the data
can be identified as the distribution function values on the spatial nodes, as well. For each cell border we compute three approximations
respectively for $\hf (\xp^{-})$ and $\hf (\xp^{+})$, as follow
\begin{equation*}
\hf^{l} (\xp^{-}) = L^{l}_{i} (\xp), \qquad \hf^{l} (\xp^{+}) = L^{l}_{i+1} (\xp), \qquad l=-1,0,1.
\end{equation*} 
Thus they are written as follow
\begin{eqnarray*}
\hf^{-1}(\xp^{-}) & = & \frac{1}{3} \brf_{i-2}  -\frac{7}{6} \brf_{i-1} \frac{11}{6} \brf_{i} \\
\hf^{0}(\xp^{-}) & = & -\frac{1}{6} \brf_{i-1} + \frac{5}{6} \brf_{i} + \frac{1}{3} \brf_{i+1} \\
\hf^{1}(\xp^{-}) & = & \frac{1}{3}\brf_{i} + \frac{5}{6} \brf_{i+1} - \frac{1}{6} \brf_{i+2} .  
\end{eqnarray*}
The computation of $\hf^{l}(\xp^{+})$ is performed by using the stencils symmetrically mirrored.  
To compute the weights we need the smoothness indicators 
\begin{eqnarray*}
\beta_i^{-1} & = & \frac{13}{12}(\brf_{i-2} -2\brf_{i-1} + \brf_{i})^2 + \frac{1}{4}(\brf_{i-2} -4\brf_{i-1} +3\brf_{i})^2,\\
\beta_i^{0} & = & \frac{13}{12}(\brf_{i-1} -2\brf_{i} + \brf_{i+1})^2 + \frac{1}{4}(\brf_{i-1} - \brf_{i+1})^2,\\
\beta_i^{1} & = & \frac{13}{12}(\brf_{i} -2\brf_{i+1} + \brf_{i+2})^2 + \frac{1}{4}(3\brf_{i} -4\brf_{i+1} +\brf_{i+2})^2
\end{eqnarray*} 
The nonlinear weights, $w^{l}_i$ and $\tilde{w}^{l}_{i}$, are computed as in \eqref{E:nonlinearw}, but now the $\alpha$ coefficients are 
evaluated as follow
\begin{equation*} 
\alpha^{l}_{i} = \frac{d_l}{(\varepsilon + \beta^{l}_{i})^2}, \qquad 
\tilde{\alpha}^{l}_{i} = \frac{\tilde{d}_l}{(\varepsilon + \beta^{l}_{i})^2}
\end{equation*}
where
\begin{equation*}
d_{-1} = \frac{1}{10}, \quad d_{0} = \frac{3}{5}, \quad d_{1} = \frac{3}{10}, \qquad \tilde{d}_{l} = d_{-l}, \qquad l=-1,0,1. 
\end{equation*}
Finally, we have
\begin{equation} \label{eq:wflux}
\hf (\xp^{-}) = \sum_{l=-1}^{1} w^{l}_{i} \hf^{l} (\xp^{-}) , \qquad \hf (\xp^{+}) = \sum_{l=-1}^{1} \tilde{w}^{l}_{i} \hf^{l} (\xp^{+}).
\end{equation}
    
\begin{figure}[htbp]
\centering
 \includegraphics[height=6cm,width=8cm,angle=0]{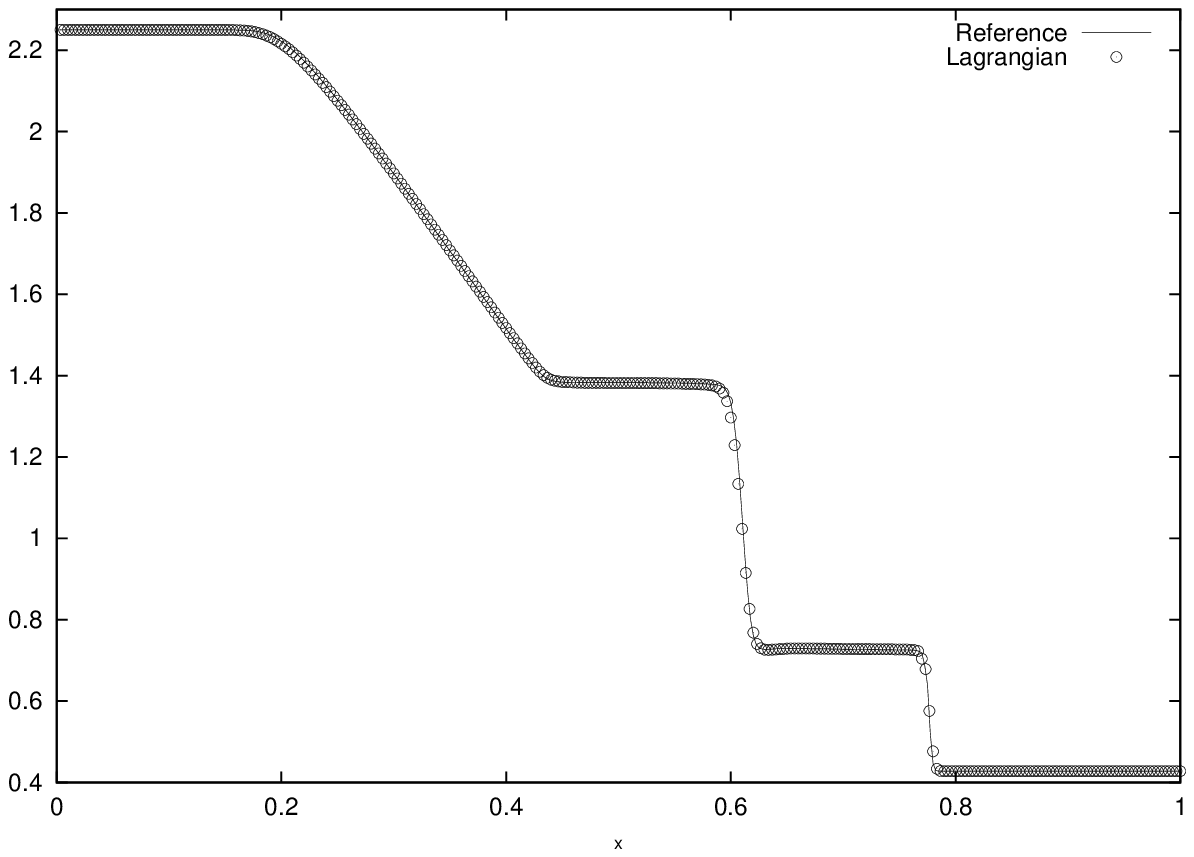}
 \includegraphics[height=6cm,width=8cm,angle=0]{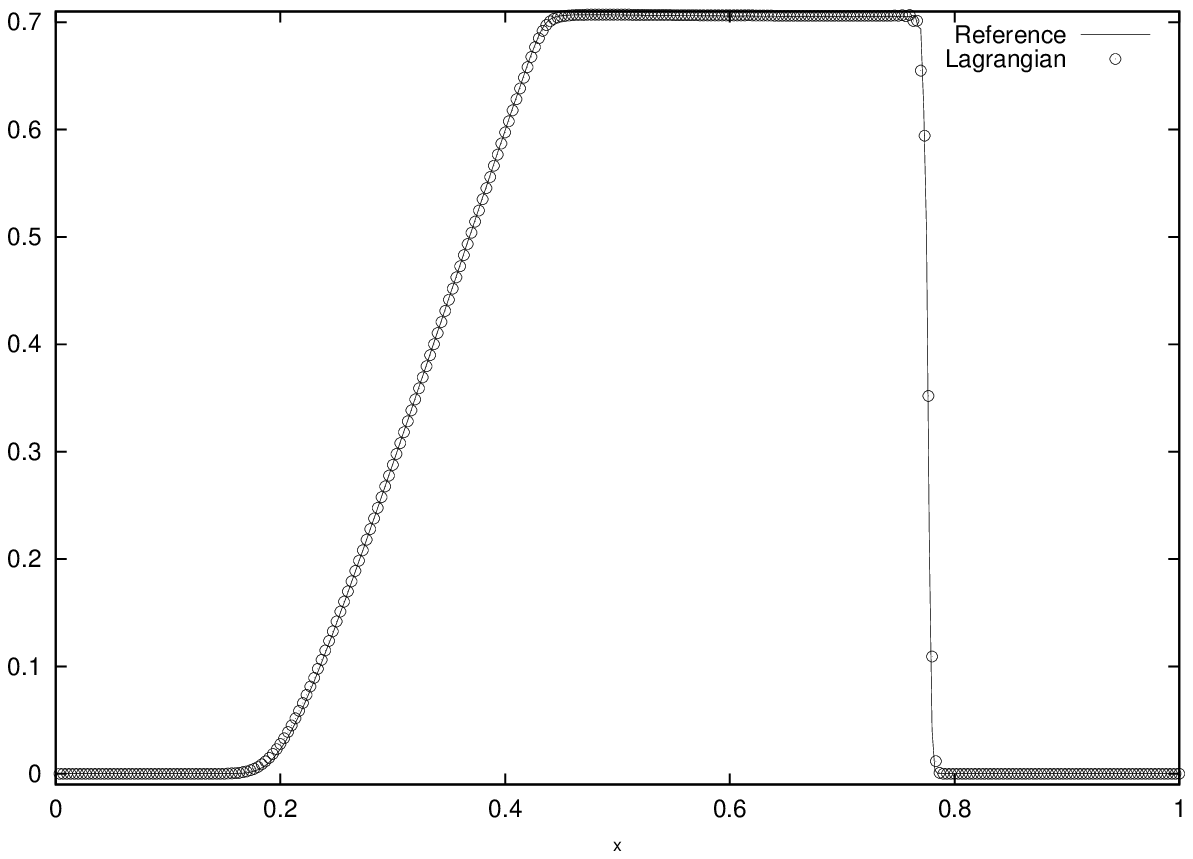}
 \includegraphics[height=6cm,width=8cm,angle=0]{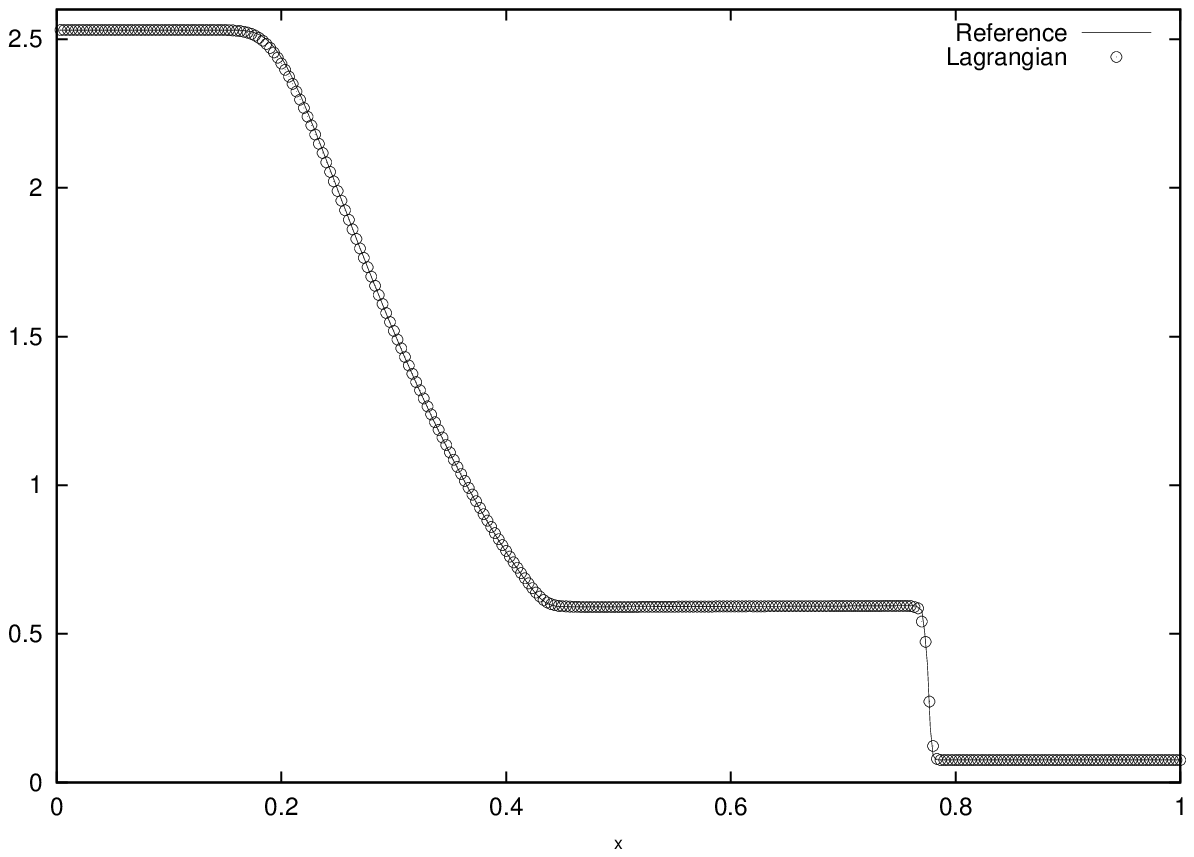}
\caption{Test 2. From the top to the bottom, density, velocity and pressure,
 $\tau=10^{-6}$. Reference solution (Nx=400) (solid line) and the 
 large time step DIRK scheme Nx=300, CFL=4.5 (star)}
 \label{fig:sod_test} 
\end{figure}

Summarizing, initially we perform a relaxation step, by using the first order implicit semi-Lagrangian scheme, to get a suitable predictor value 
of the internal stage, $F^{(1)}(v)$. This is used to compute $\hf^{l} (\xp^{-})$ and $\hf^{l} (\xp^{+})$ for each cell to get $\hf (\xp^{-})$ and 
$\hf (\xp^{+})$, by \eqref{eq:wflux} respectively.   
In this part of the computation the CFL can be greater than one, of course, 
since we know in advance the solution of the characteristic lines. Once the flux is evaluated at each cell border by \eqref{E:fluxes}, 
it is possible to compute the discrete convective term by \eqref{shuDer} . Finally the corrector step is applied by solving implicitly 
the equation \eqref{E:consForm}. This procedure can be straightforwardly extended to multiple internal stages. 
\figurename~\ref{fig:sod_test} shows the comparison between the proposed numerical method and the reference solution computed by a 
high order solver of the Euler equations. We observe a clear improvement of the quality of the solution. 
Also in this case the accuracy, despite a slight improvement, doesn't hold the third order in fluid regime, but the conservation 
is fully achieved.    

\begin{table}[!ht] 
\begin{center}
\begin{tabular}{|r|c|c|c|}
\cline{1-4}
\multicolumn{4}{|c|}{$L_{2}-Relative \, errors$} \\ \hline
 $N_x$ & Density & Velocity  & Temperature \\ \hline
 20 & 2.54838e-03 & 2.35049e-03 & 4.63423e-03 \\
 40 & 5.57339e-04 & 3.64146e-04 & 9.17053e-04 \\
 80 & 8.41532e-05 & 5.21314e-05 & 1.28062e-04 \\
160 & 1.17817e-05 & 8.94658e-06 & 2.53109e-05 \\
320 & 1.69746e-06 & 1.95126e-06 & 6.47027e-06 \\ 
 \hline
\multicolumn{4}{|c|}{$L_{2}-Orders$} \\ \hline
$N_x$ & Density & Velocity  & Temperature \\ \hline
 40 & 2.193 & 2.690 & 2.337 \\
 80 & 2.727 & 2.804 & 2.840 \\
160 & 2.836 & 2.543 & 2.339 \\
320 & 2.795 & 2.197 & 1.968 \\
 \hline
\end{tabular}
\end{center}
\caption{Scheme S2, $\tau=10^{-2},CFL=4.5$.} 
\label{tab:err1}
\end{table}

\begin{table}[!ht] 
\begin{center}
\begin{tabular}{|r|c|c|c|}
\cline{1-4}
\multicolumn{4}{|c|}{$L_{2}-Relative \, errors$} \\ \hline
 $N_x$ & Density & Velocity  & Temperature \\ \hline
 20 & 2.96809e-03 & 3.15227e-03 & 6.37101e-03 \\
 40 & 6.57722e-04 & 6.05216e-04 & 1.94043e-03 \\ 
 80 & 1.11120e-04 & 1.36059e-04 & 5.26168e-04 \\
160 & 2.25137e-05 & 4.61239e-05 & 1.59816e-04 \\
320 & 6.10643e-06 & 1.63245e-05 & 5.05549e-05 \\
 \hline
\multicolumn{4}{|c|}{$L_{2}-Orders$} \\ \hline
$N_x$ & Density & Velocity  & Temperature \\ \hline
 40 & 2.174 & 2.381 & 1.715 \\
 80 & 2.565 & 2.153 & 1.883 \\
160 & 2.303 & 1.561 & 1.719 \\
320 & 1.882 & 1.498 & 1.660 \\
 \hline
\end{tabular}
\end{center}
\caption{Scheme S2,$\tau=10^{-6},CFL=4.5$.} 
\label{tab:err2}
\end{table}
 
\begin{table}[!ht]
\begin{center}
\begin{tabular}{|r|c|c|c|}
\cline{1-4}
\multicolumn{4}{|c|}{$L_{2}-Relative \, errors$} \\ \hline
 $N_x$ & Density & Velocity  & Temperature \\ \hline
 20 & 2.41539e-03 & 1.97185e-03 & 4.36445e-03 \\ 
 40 & 4.93444e-04 & 2.90747e-04 & 8.14397e-04 \\   
 80 & 7.36995e-05 & 4.27296e-05 & 1.14397e-04 \\
160 & 1.06248e-05 & 5.91413e-06 & 1.54660e-05 \\ 
320 & 1.55051e-06 & 9.55269e-07 & 2.89414e-06 \\ 
 \hline
\multicolumn{4}{|c|}{$L_{2}-Orders$} \\ \hline
$N_x$ & Density & Velocity  & Temperature \\ \hline
 40 & 2.291 & 2.762 & 2.422 \\
 80 & 2.743 & 2.766 & 2.832 \\
160 & 2.794 & 2.853 & 2.887 \\
320 & 2.777 & 2.630 & 2.418 \\
 \hline
\end{tabular}
\end{center}
\caption{Scheme S3,$\tau=10^{-2},CFL=4.5$.} 
\label{tab:err3} 
\end{table}
\begin{table}[!ht]
\begin{center}
\begin{tabular}{|r|c|c|c|}
\cline{1-4}
\multicolumn{4}{|c|}{$L_{2}-Relative \, errors$} \\ \hline
 $N_x$ & Density & Velocity  & Temperature \\ \hline
 20 & 2.02024e-03 & 2.72023e-03 & 4.79517e-03 \\
 40 & 5.04007e-04 & 7.57946e-04 & 2.06197e-03 \\
 80 & 9.91878e-05 & 2.63921e-04 & 9.43964e-04 \\
160 & 4.10972e-05 & 1.48278e-04 & 5.14779e-04 \\
320 & 2.17256e-05 & 7.52320e-05 & 2.39514e-04 \\
 \hline
\multicolumn{4}{|c|}{$L_{2}-Orders$} \\ \hline
$N_x$ & Density & Velocity  & Temperature \\ \hline
 40 & 2.003 & 1.844 & 1.218 \\
 80 & 2.345 & 1.522 & 1.127 \\
160 & 1.271 & 0.932 & 0.975 \\
320 & 1.120 & 1.079 & 1.104 \\
 \hline
\end{tabular}
\end{center}
\caption{Scheme S3,$\tau=10^{-6},CFL=4.5$.} 
\label{tab:err4}
\end{table}
 
\begin{table}[!ht]
\begin{center}
\begin{tabular}{|r|c|c|c|}
\cline{1-4}
 \multicolumn{4}{|c|}{$CFL = 1.5$}  \\ \hline  
    & Density & Velocity  & Temperature \\  \hline 
 S1 & 6.000e-05 & 1.120e-04 & 1.700e-04 \\
 S2 & 6.000e-05 & 1.120e-04 & 1.700e-04 \\
 S3 & 6.000e-05 & 1.120e-04 & 1.700e-04 \\ 
 \hline 
 \multicolumn{4}{|c|}{$CFL = 5.5$} \\ \hline 
    & Density & Velocity  & Temperature \\ \hline
 S1 & 2.500e-04 & 2.977e-04 & 4.400e-04 \\
 S2 & 2.500e-04 & 2.972e-04 & 4.400e-04 \\
 S3 & 2.500e-04 & 2.973e-04 & 4.400e-04 \\ 
 \hline 
 \multicolumn{4}{|c|}{$CFL = 10.5$} \\ \hline
    & Density & Velocity  & Temperature \\ \hline
 S1 & 3.000e-04 & 3.120e-04 & 4.700e-04 \\
 S2 & 3.000e-04 & 3.115e-04 & 4.700e-04 \\
 S3 & 3.000e-04 & 3.117e-04 & 4.700e-04 \\ 
 \hline
\end{tabular}
\end{center}
\caption{Comparison test for different CFL. Test 1; $N_x = 200$, $\tau = 10^{-2}$.} 
\label{tab:cflComp}
\end{table}

\begin{table}[!ht]
\begin{center}
\begin{tabular}{|r|c|c|c|}
\cline{1-4}
\multicolumn{4}{|c|}{$\tau=$1e-2} \\ \hline
$N_x$ & Density & Momentum  & Energy \\ \hline
  20 & 3.19103e-04 & 6.45517e-04 & 6.44489e-04 \\
  40 & 8.68364e-05 & 2.43064e-05 & 1.53479e-04 \\ 
  80 & 1.86619e-05 & 6.93736e-06 & 3.38482e-05 \\
 160 & 2.21369e-06 & 6.03659e-07 & 3.82991e-06 \\
 320 & 2.47089e-07 & 6.48153e-08 & 4.23732e-07 \\
 640 & 2.75503e-08 & 5.92434e-09 & 4.57648e-08 \\
1280 & 3.21756e-09 & 6.82768e-10 & 5.29436e-09 \\ 
\hline
\cline{2-3}
\multicolumn{4}{|c|}{$\tau=$1e-6} \\ \hline
$N_x$ & Density & Momentum  & Energy \\ \hline
  20 & 3.86571e-04 & 8.49545e-04 & 8.59503e-04 \\
  40 & 1.25170e-04 & 4.22174e-05 & 2.34721e-04 \\
  80 & 3.17682e-05 & 1.53056e-05 & 6.50680e-05 \\
 160 & 4.65888e-06 & 1.86177e-06 & 9.55101e-06 \\
 320 & 5.94587e-07 & 2.25448e-07 & 1.20946e-06 \\
 640 & 7.47884e-08 & 2.69514e-08 & 1.52337e-07 \\
1280 & 9.24536e-09 & 3.27055e-09 & 1.87483e-08 \\
\hline
\end{tabular}
\end{center}
\caption{Errors in conservation. Test 1, scheme S2.} 
\label{tab:conserv}
\end{table}
\newpage 
\section{Conclusions}
In this work a new class of semi-Lagrangian schemes for BGK model of the Boltzmann equation has been introduced. 
These schemes show to be highly efficient and accurate, allowing us, not only to investigate the rarefied regime
of the gasdynamics system, but also to converge correctly to the fluid limit for discontinous solution (shock wave propagation).
Moreover also high CFL values (less than one in classical schemes) do not affect the numerical solution by the diffusivity, preserving 
accuracy and conservation. These properties make the extension to multidimensional schemes highly desirable, leading up to
the possibility to simulate more realistic problems, like bidimensional generalized Riemann problems or microfluidics-based devices 
of engineering interest. The authors are currently go along this way, and a conservative multidimensional version 
of this class of schemes is under investigation.   
\newpage 
\section*{Ackowledgemts}
The authors would like to thank the support received 
for project MONUMENT (MOdellizzazione NUmerica in MEms and NanoTecnologie) by
Galileo Programme of the \emph{Universit\`a Italo Francese},
institution engaged for the academic cooperation between
Italy and France.
 
\begin{figure}[b]
\centering
 \includegraphics[height=6cm,width=6.5cm,angle=0]{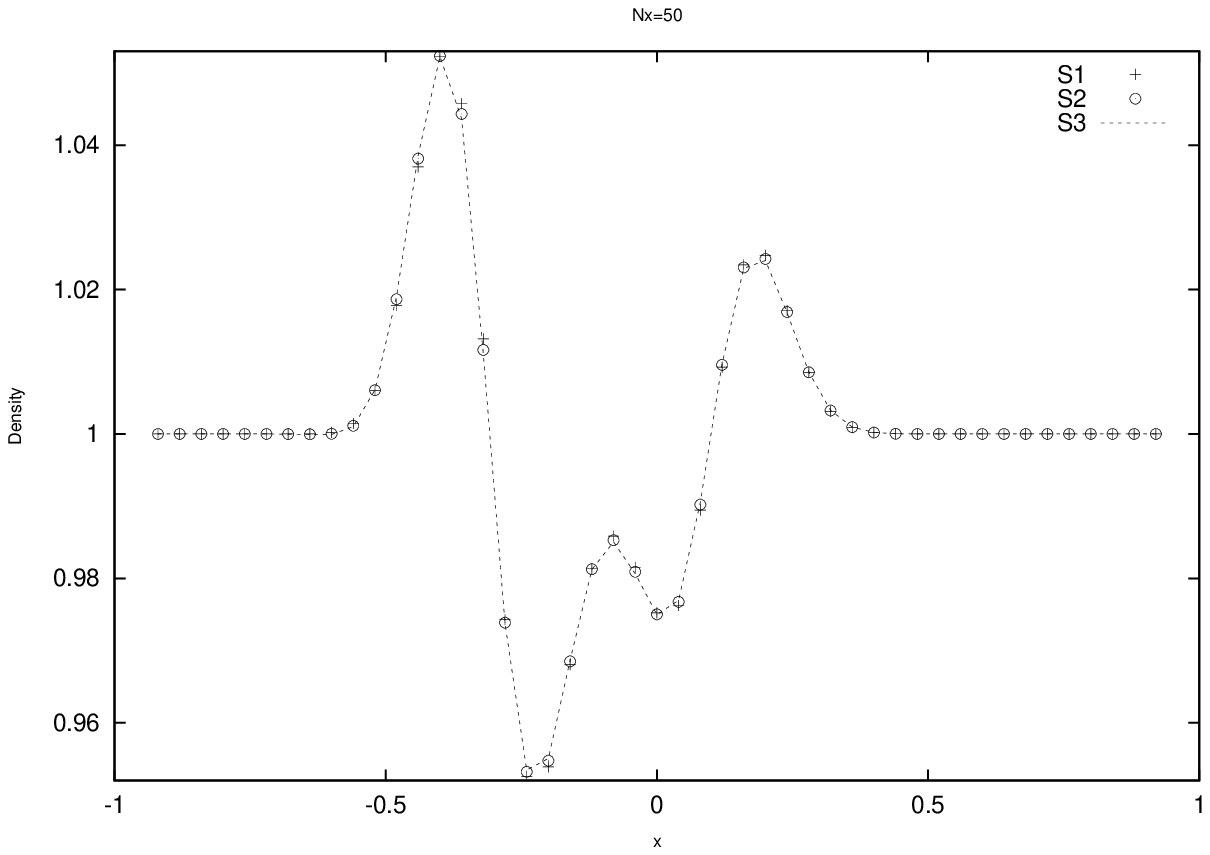}%
 \includegraphics[height=6cm,width=6cm,angle=0]{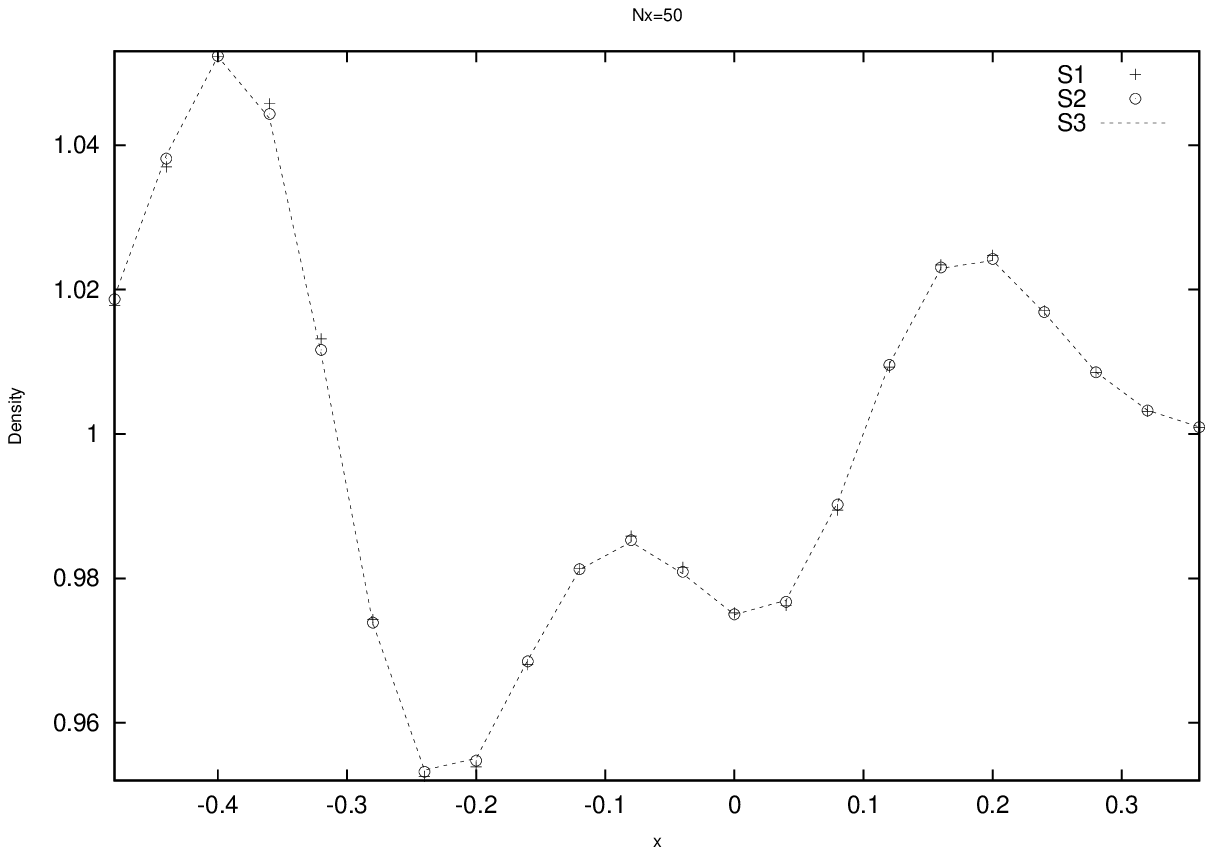}
 \includegraphics[height=6cm,width=6.5cm,angle=0]{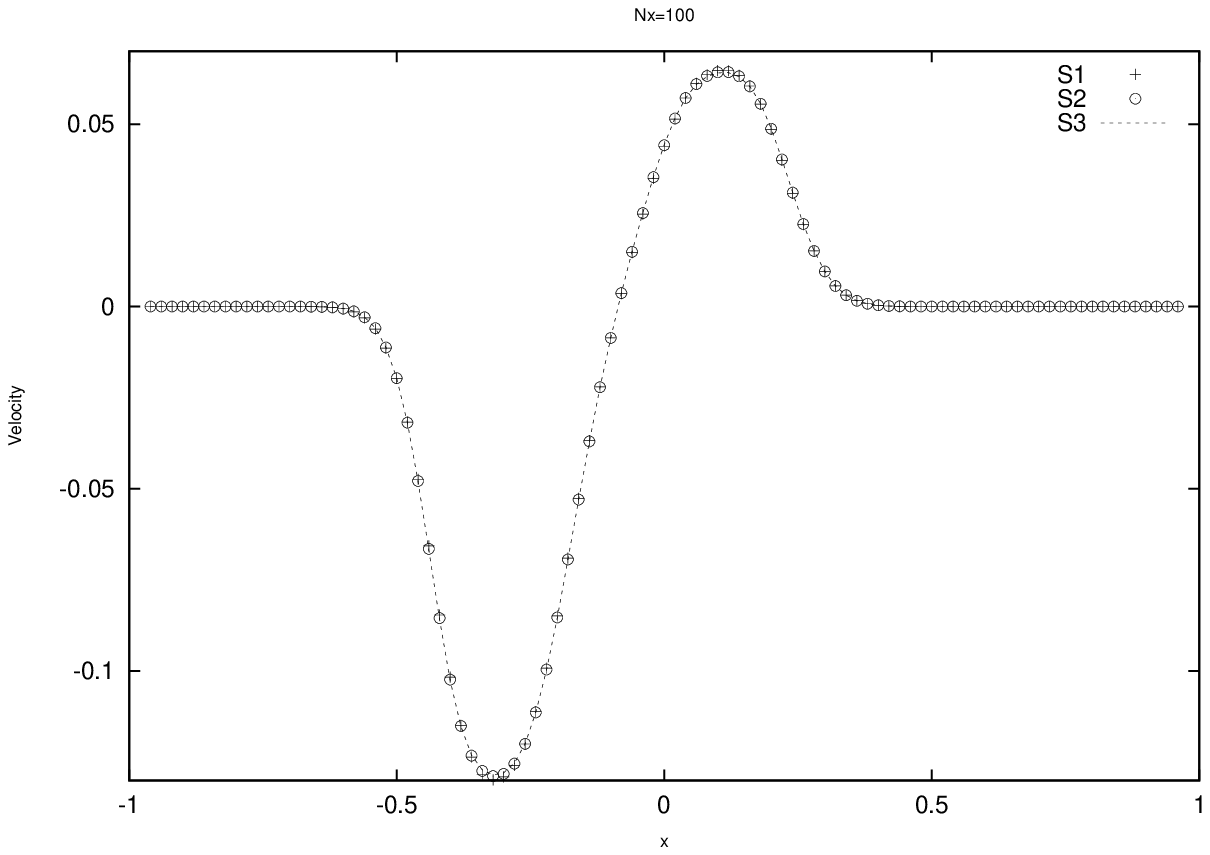}
 \includegraphics[height=6cm,width=6cm,angle=0]{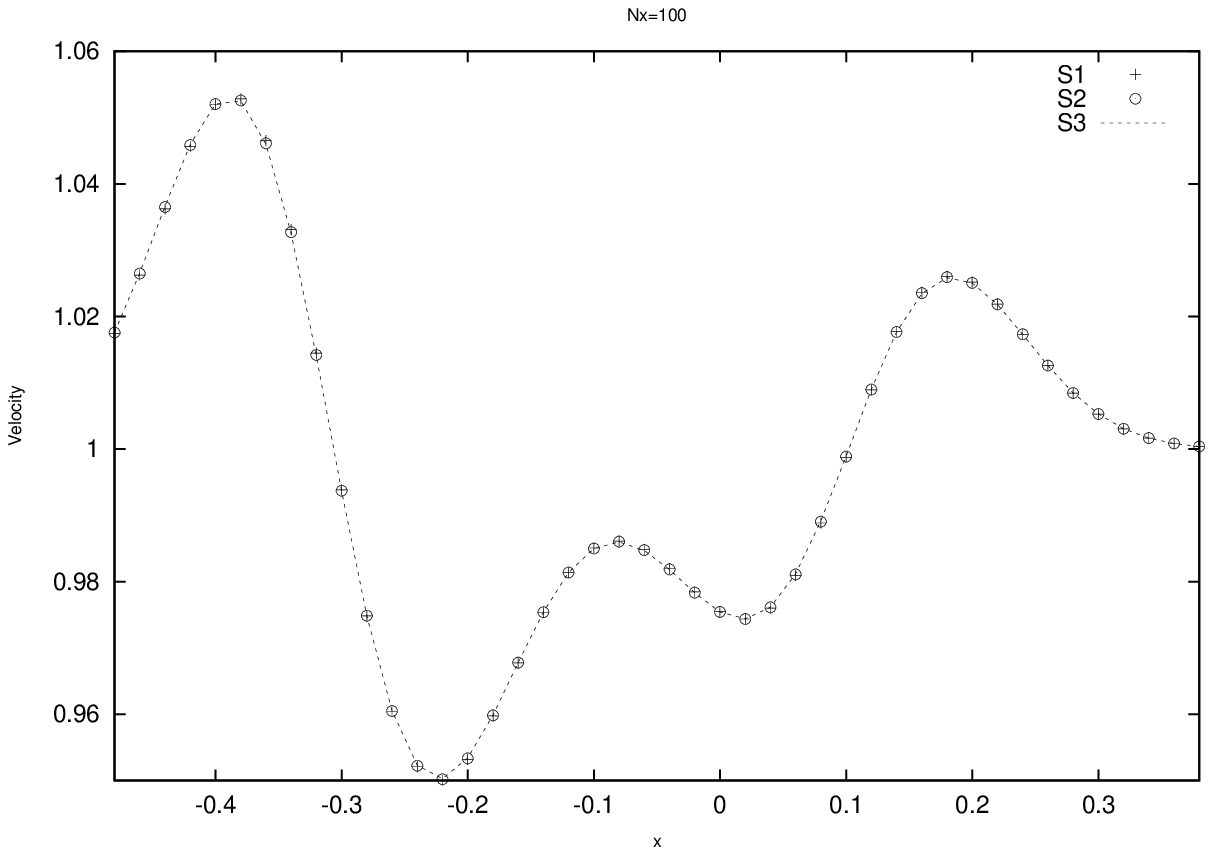}
 \includegraphics[height=6cm,width=6.5cm,angle=0]{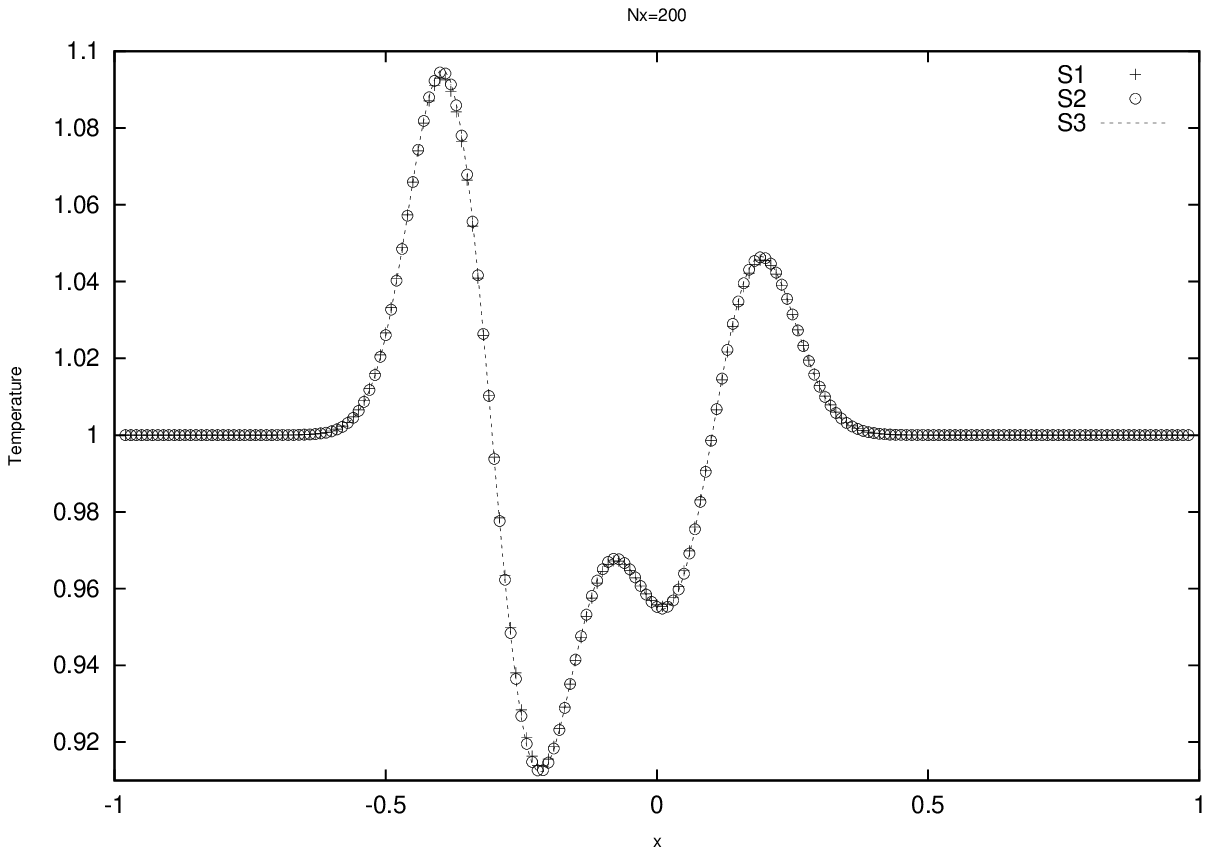}
 \includegraphics[height=6cm,width=6cm,angle=0]{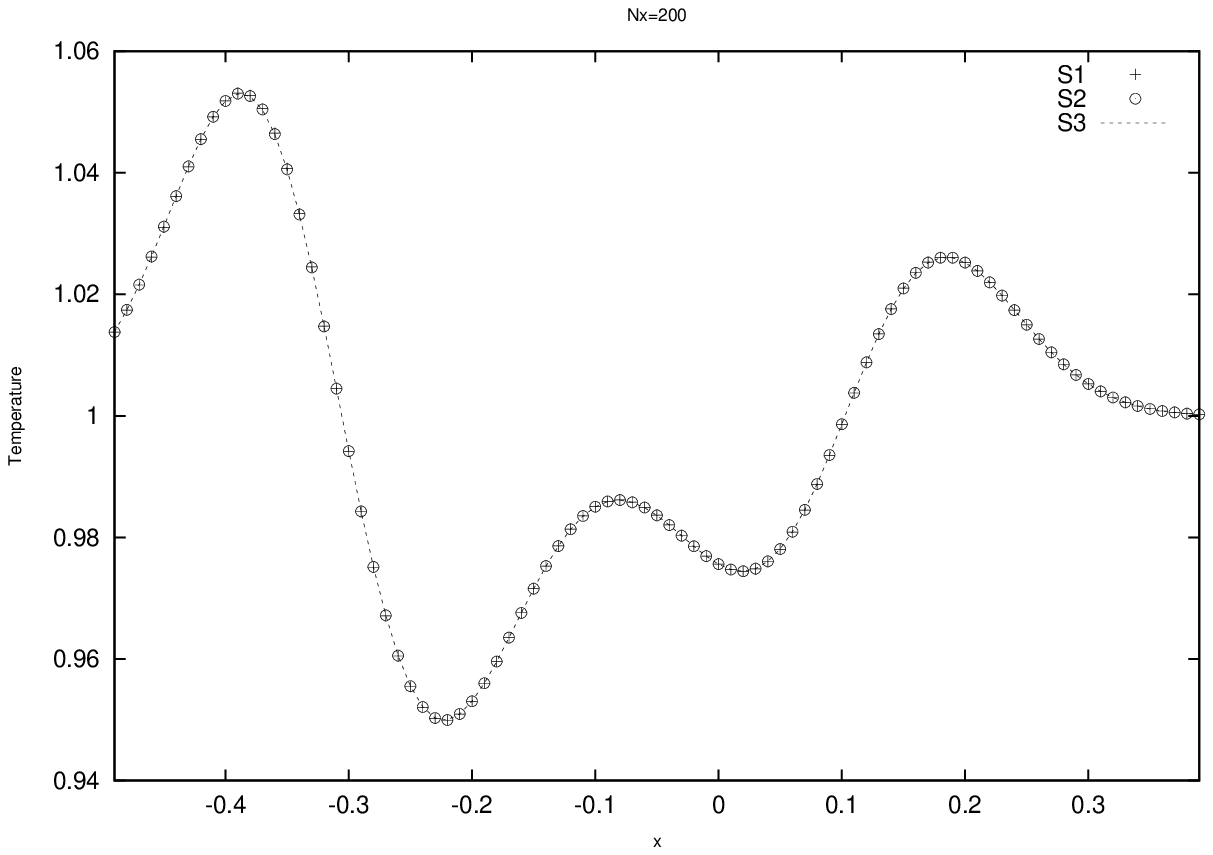}
 \caption{Test 1,$\tau=10^{-6}$,CFL=10.0 . Left column: from the top 
          to the bottom, density, velocity and temperature for the schemes S1 (cross),
          S2 (circle), S3 (dashed). Right column: from the top to the bottom, density zoom 
          for $N_x=40,80,160$, respectively.}
\end{figure}
\begin{figure}[htbp]
\centering
 \includegraphics[height=6cm,width=6.25cm,angle=0]{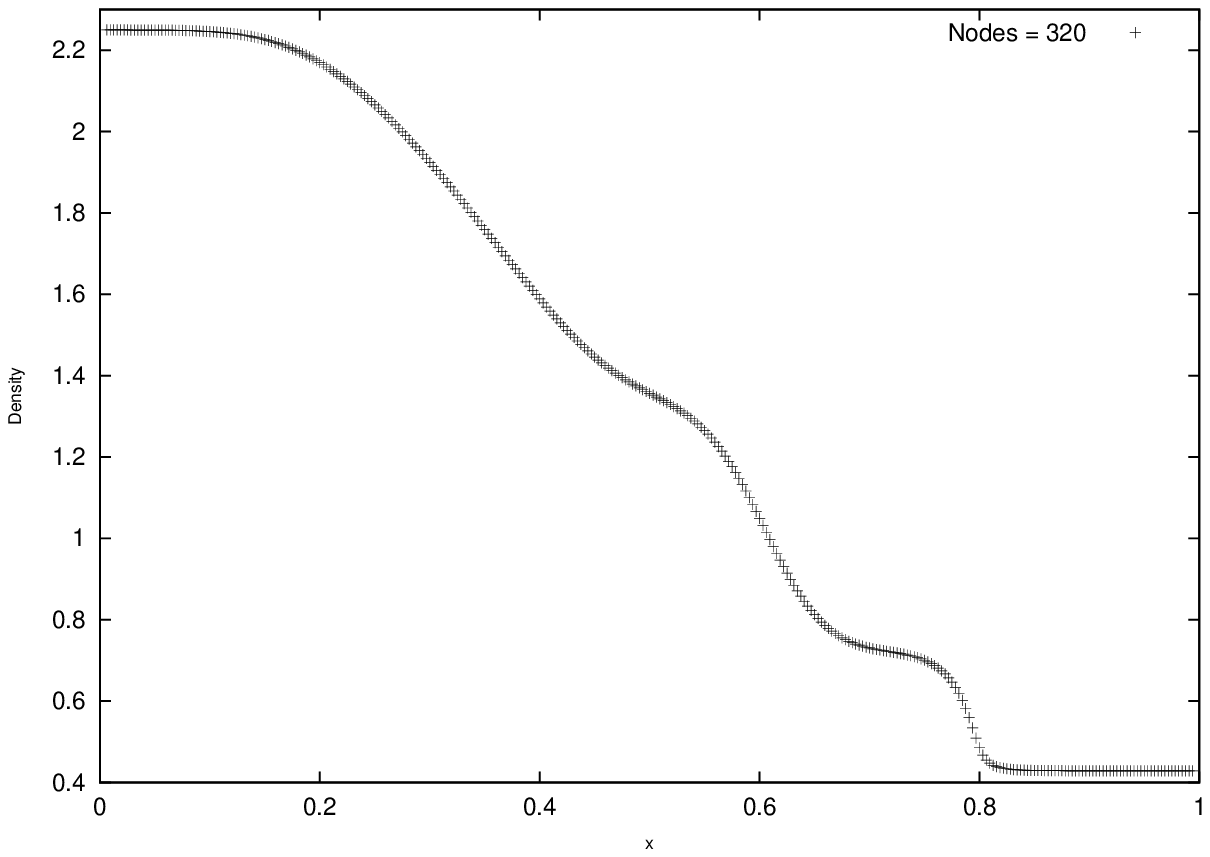}%
 \includegraphics[height=6cm,width=6.25cm,angle=0]{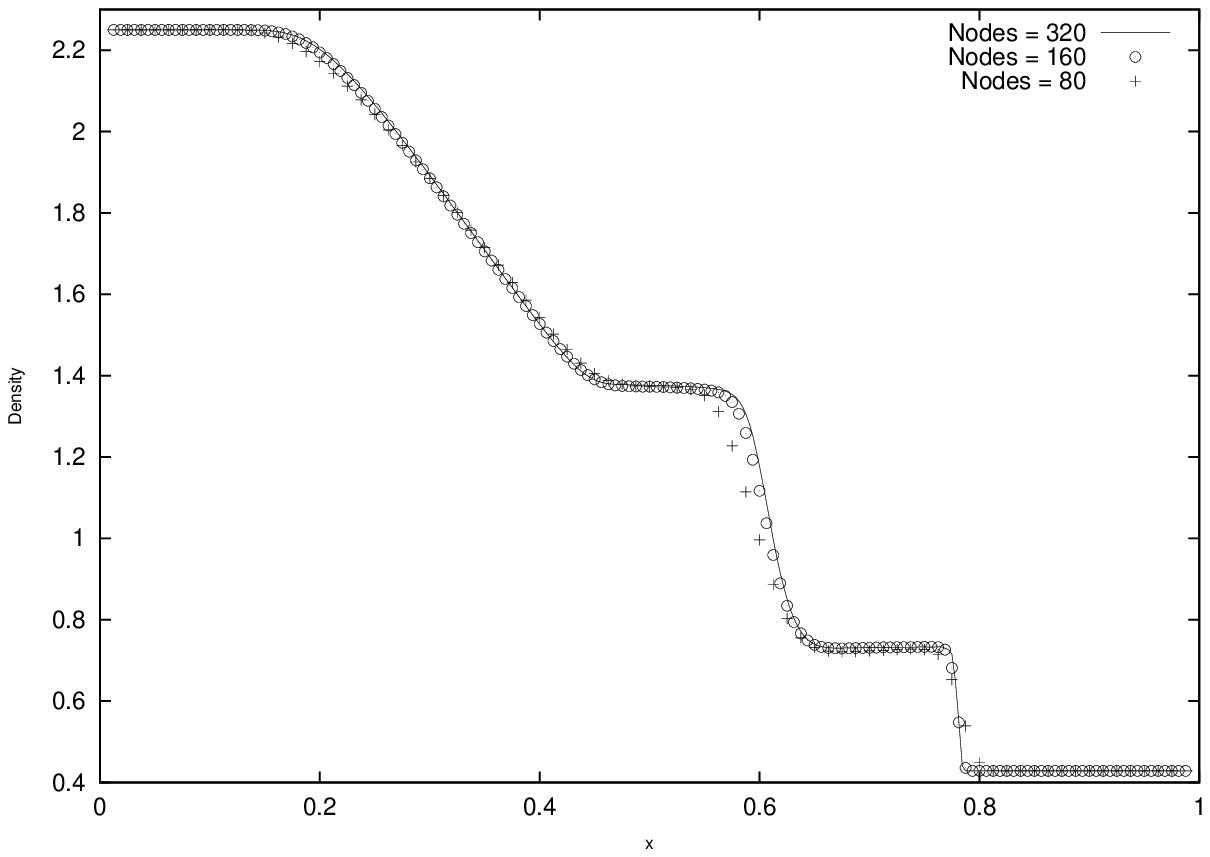}
 \includegraphics[height=6cm,width=6.25cm,angle=0]{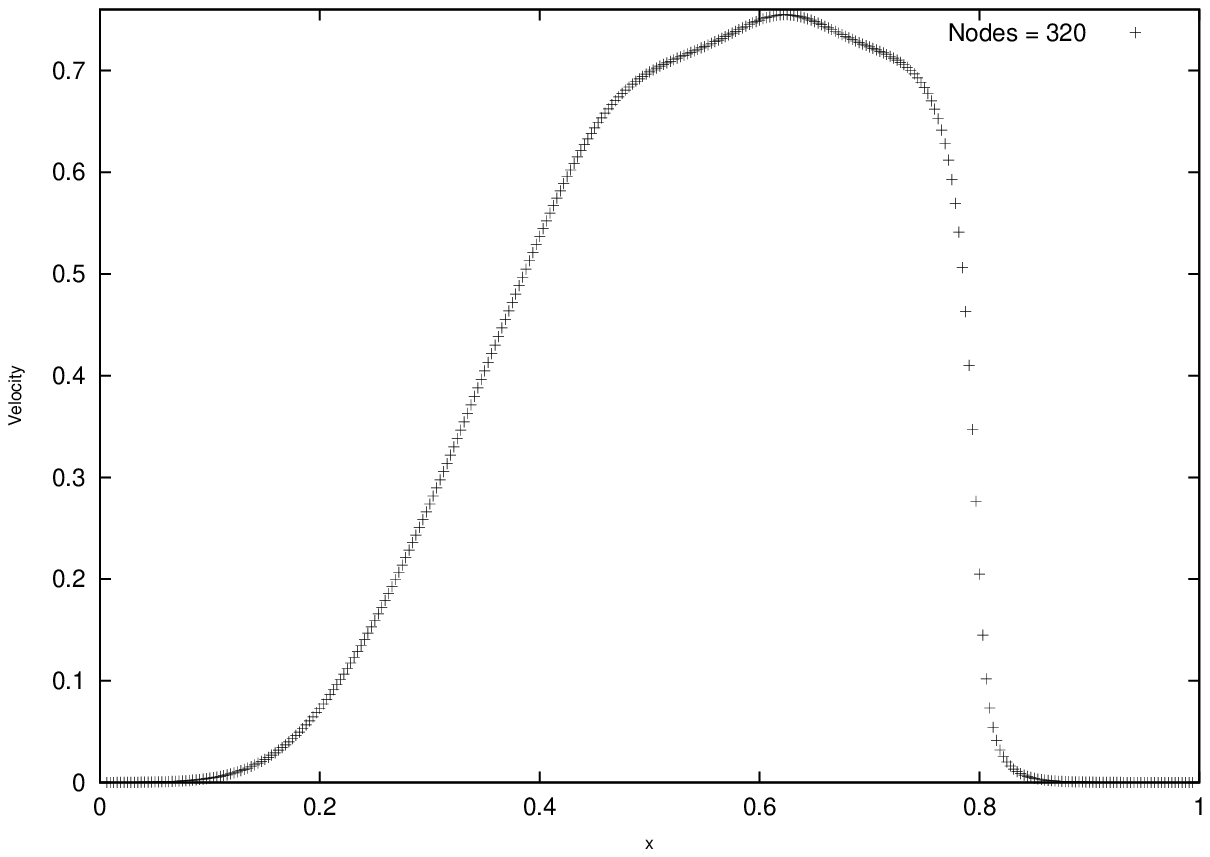}
 \includegraphics[height=6cm,width=6.25cm,angle=0]{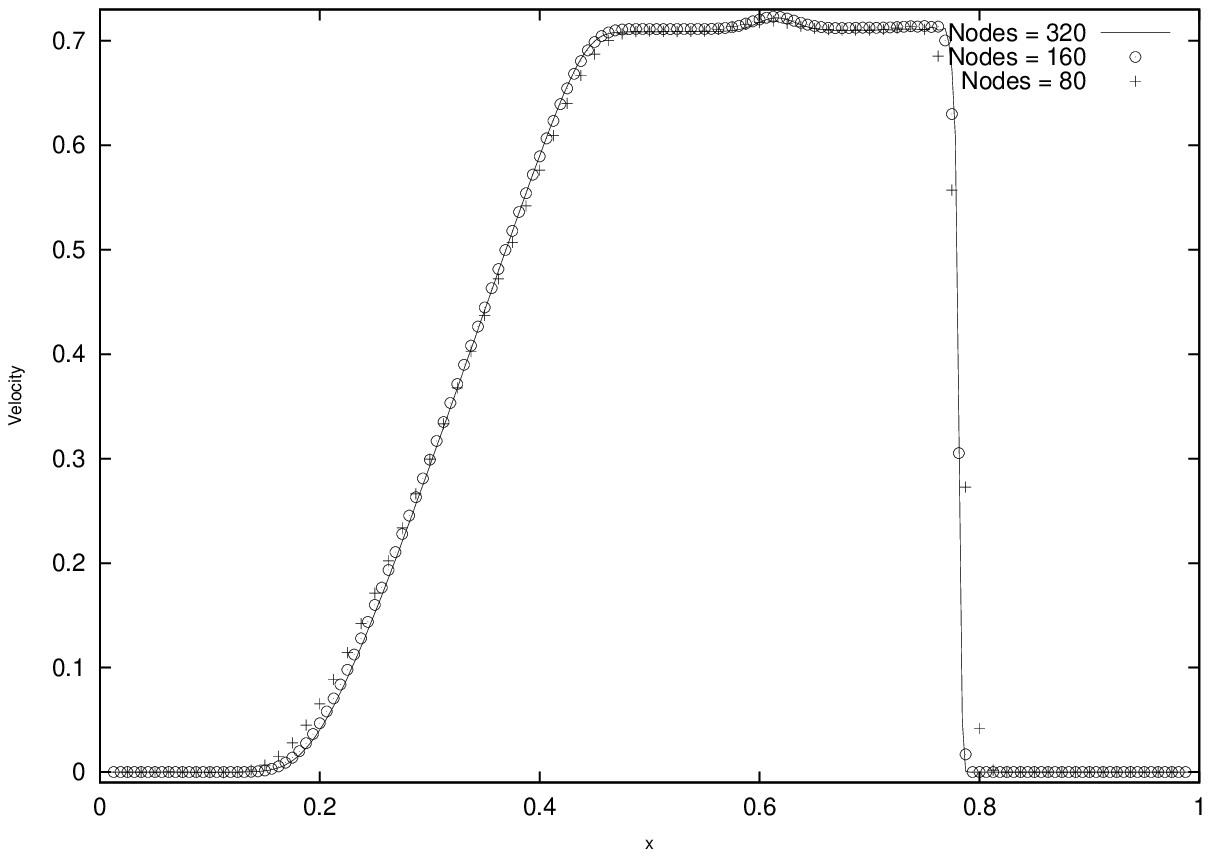}
 \includegraphics[height=6cm,width=6.25cm,angle=0]{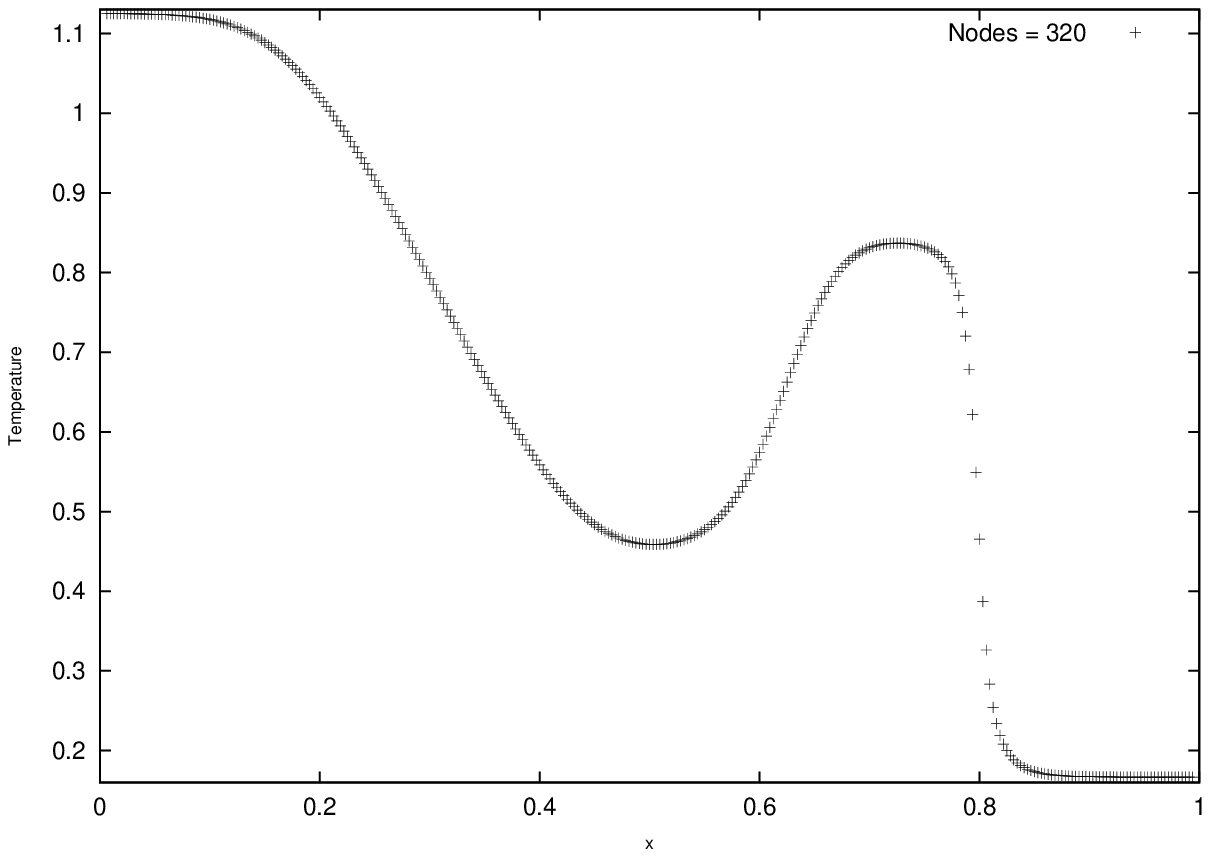}
 \includegraphics[height=6cm,width=6.25cm,angle=0]{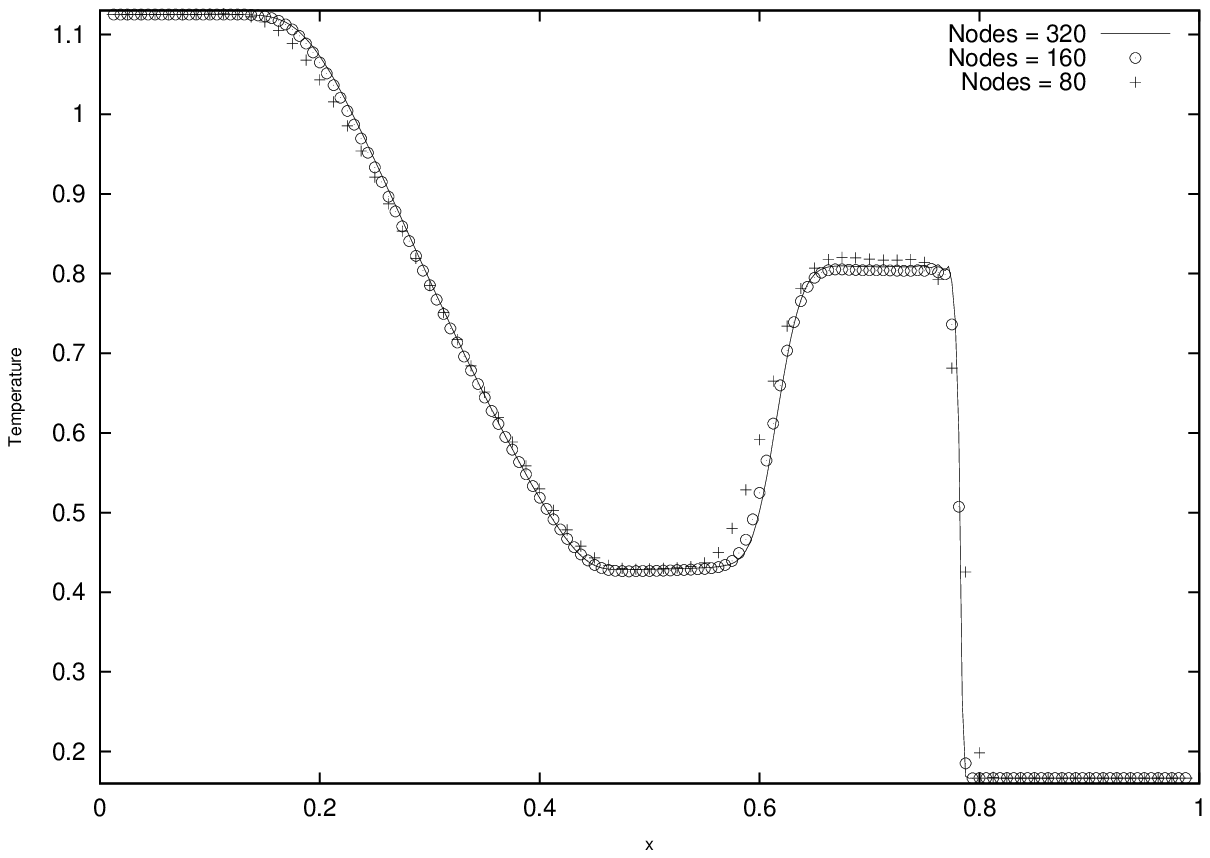} 
\caption{Test 2. From the top to the bottom, density, velocity and temperature. Left column
 $\tau=10^{-2}$. Right column $\tau=10^{-6}$. CFL=9.44.}
\end{figure}


\end{document}